\documentclass[aos,citesort,dvips]{arximspdf}
\usepackage{mathrsfs,graphics}

\doi{10.1214/10-AOS799}
\volume{38}
\issue{5}
\pubyear{2010}
\firstpage{2916}
\lastpage{2957}

\makeatletter

\newcommand{\iint}{\int\!\!\int}

\newproclaim{remark}{Remark}
\newproclaim{example}{Example}
\newtheorem{lemma}{Lemma}
\newtheorem{proposition}{Proposition}
\newtheorem{theorem}{Theorem}
\newproclaim{algorithm}{Algorithm}

\newcommand{\e}{{e}}

\makeatother

\begin{document}
\begin{frontmatter}

\title{Kernel density estimation via diffusion}
\runtitle{Kernel density estimation via diffusion}

\begin{aug}
\author[A]{\fnms{Z. I.} \snm{Botev}\corref{}\thanksref{t1}\ead[label=e1]{botev@maths.uq.edu.au}},
\author[A]{\fnms{J. F.} \snm{Grotowski}\ead[label=e2]{grotow@maths.uq.edu.au}} and
\author[A]{\fnms{D. P.} \snm{Kroese}\thanksref{t1}\ead[label=e3]{kroese@maths.uq.edu.au}\ead[label=u1,url]{http://www.maths.uq.edu.au/\textasciitilde kroese/}}
\runauthor{Z. I. Botev, J. F. Grotowski and D. P. Kroese}
\affiliation{University of Queensland}
\address[A]{School of Mathematics and Physics\\
University of Queensland\\
St. Lucia, Brisbane\\
Queensland, 4072\\
Australia\\
\printead{e1}\\
\phantom{E-mail: }\printead*{e2}\\
\phantom{E-mail: }\printead*{e3}\\
\printead{u1}} 
\end{aug}

\thankstext{t1}{Supported by Australian Research
Council Grant DP0985177.}

\received{\smonth{12} \syear{2009}}

%
\begin{abstract}
We present a new adaptive kernel density estimator based on linear
diffusion processes. The proposed estimator builds on existing ideas
for adaptive smoothing by incorporating information from a pilot
density estimate. In addition, we propose a new plug-in bandwidth
selection method that is free from the arbitrary normal reference rules
used by existing methods. We present simulation examples in which
the proposed approach outperforms existing methods in terms of
accuracy and reliability.
\end{abstract}

%

\end{frontmatter}

\section{Introduction}\label{sec1}
Nonparametric density estimation is an important tool in the
statistical analysis of data. A nonparametric estimate can be used, for
example, to assess the
multimodality, skewness, or any other structure in the distribution of
the data \cite{Scott,Silverman}.
It can also be used
for the summarization of Bayesian posteriors, classification and
discriminant analysis \cite{Simonoff}. Nonparametric density
estimation has even proved useful in Monte Carlo computational methods,
such as the smoothed bootstrap method
and the particle filter method \cite{doufregor01}. Nonparametric
density estimation is an alternative to the
parametric approach, in which one specifies a model
up to a small number of parameters and then estimates the parameters
via the likelihood principle. The
advantage of the nonparametric approach is that it offers a
far greater flexibility in modeling a given dataset and, unlike the
classical approach, is not affected by specification bias
\cite{lehmann90}. Currently, the most popular nonparametric approach to density
estimation is \textit{kernel density estimation} (see
\cite{Scott,Simonoff,WandJones}).

Despite the vast body of literature on the subject,
there are still many contentious issues regarding the
implementation and practical performance of kernel density
estimators.
First, the most popular data-driven bandwidth selection technique, the
\textit{plug-in} method \cite{sheather96,sheather91}, is adversely
affected by the so-called
\textit{normal reference rule} \cite{devroye97,jones96}, which is
essentially a construction of a preliminary normal
model of the data upon which the performance of the bandwidth
selection method depends. Although plug-in estimators perform
well when the normality assumption holds approximately,
at a conceptual level the use of the normal reference rule invalidates the
original motivation for applying
a nonparametric method in the first place.

Second, the popular Gaussian kernel density estimator \cite{MarronWand}
lacks
\textit{local adaptivity}, and this often results in a large sensitivity
to outliers, the presence of spurious bumps, and in an
overall unsatisfactory bias performance---a tendency to flatten the
peaks and valleys of the density \cite{TerrelScott}.

Third, most kernel estimators suffer from \textit{boundary bias} when,
for example, the
data is nonnegative---a phenomenon due to the fact that most kernels
do not take into account specific knowledge about the domain of the
data \cite{marron94,park03}.

These problems have been alleviated to a certain degree by the
introduction of more sophisticated kernels than the simple Gaussian kernel.
Higher-order kernels have been used as a way to improve local
adaptivity and reduce bias \cite{jones97}, but these have the
disadvantages of not giving
proper nonnegative density estimates, and of requiring a large sample
size for good performance \cite{MarronWand}.
The lack of local adaptivity has been addressed
by the introduction of \textit{adaptive} kernel estimators
\cite{abramson,hall90,hall95,jones94}. These include the
\textit{balloon} estimators, \textit{nearest
neighbor} estimators and \textit{variable bandwidth} kernel estimators
\cite{loftsgaarden,TerrelScott}, none of which yield bona fide densities,
and thus remain somewhat unsatisfactory. Other
proposals such as the \textit{sample point adaptive} estimators are
computationally burdensome (the fast Fourier
transform cannot be applied \cite{Silverman}),
and in some cases do not integrate to unity \cite{park03}. The
\textit{boundary kernel estimators} \cite{jones93}, which are specifically
designed to deal with boundary bias, are either not adaptive away from
the boundaries or do not result in bona fide
densities \cite{foster}. Thus, the literature abounds with partial
solutions that obscure a unified comprehensive framework for the
resolution of these problems.\looseness=1

The aim of this paper is to introduce an
adaptive kernel density estimation method
based on the smoothing properties of linear diffusion processes. The
key idea is to view the kernel from which the estimator is constructed
as the transition density of a diffusion process.
We utilize the most
general linear diffusion process that has a
given limiting and stationary probability density. This
stationary density is selected to be either a pilot
density estimate or a density that the
statistician believes
represents the information about the data
prior to observing the available empirical
data. The approach leads to a simple and intuitive
kernel estimator with substantially reduced
asymptotic bias and mean square error. The
proposed estimator deals well with boundary bias and, unlike other
proposals, is always a bona fide probability density function.
We show that the proposed approach
brings under a single framework some well-known bias reduction methods,
such as the~Abramson estimator \cite{abramson}
and other variable location or scale estimators
\cite{choi,hall02,Samiuddin,jones94}.

In addition, the paper introduces an improved
plug-in bandwidth selection method that completely avoids
the normal reference rules \cite{jones96} that have adversely affected the
performance of plug-in methods. The new plug-in method is thus genuinely
``nonparametric,'' since it does not require a preliminary normal model
for the data. Moreover,
our plug-in approach does not involve numerical optimization and is not
much slower than computing a normal reference rule \cite{Matlab}.

The rest of the paper is organized as follows.
First, we describe the Gaussian kernel density
estimator and explain how it can be viewed as a
special case of smoothing using a diffusion
process. The Gaussian kernel density estimator is then
used to motivate the most general linear
diffusion that will have a set of essential smoothing
properties. We analyze the
asymptotic properties of the resulting estimator and explain how to
compute the asymptotically optimal plug-in bandwidth.
Finally, the practical benefits of the model are
demonstrated through simulation examples on some well-known
datasets \cite{MarronWand}. Our findings demonstrate an improved bias
performance and low computational cost, and a boundary bias improvement.

\section{Background}\label{sec2}
Given $N$ independent realizations $\mathcal{X}_N\equiv\{ X_1,\ldots
, X_N\}$ from an unknown continuous probability density function (p.d.f.)
$f$ on $\mathscr{X}$, the \textit{Gaussian kernel density estimator} is
defined as
%
%
\begin{equation}\label{KDE}
\hat f(x;t)=\frac{1}{N}\sum_{i=1}^N \phi(x,X_i;t),\qquad x\in\mathbb{R},
\end{equation}
where
\[
\phi(x,X_i;t)=\frac{1}{\sqrt{2\pi t}} \e^{-{(x-X_i)^2}/({2t})}
\]
is a Gaussian p.d.f. (kernel) with location $X_i$ and scale $\sqrt{t}$.
The scale is usually referred to as the \textit{bandwidth}. Much
research has been focused on the optimal choice of $t$
in (\ref{KDE}), because the performance of $\hat f$ as an estimator of $f$
depends crucially on its value \cite{sheather96,sheather91}. A well-studied
criterion used to determine an optimal $t$ is the \textit{Mean Integrated
Squared Error} (MISE),
\[
\operatorname{MISE}\{\hat f\}(t)=\mathbb E_f \int[\hat f(x;t)-f(x)]^2 \,dx,
\]
which is conveniently decomposed into integrated squared bias and
integrated variance components:
\[
\operatorname{MISE}\{\hat f\}(t)=
\int\bigl(\underbrace{\mathbb E_f[\hat f(x;t)]-f(x) }_{\mathrm{pointwise}\
\mathrm{bias}\ \mathrm{of}\ f }\bigr)^2 \,dx+ \int\underbrace{\operatorname
{Var}_f[\hat f(x;t)]
}_{\mathrm{pointwise}\ \mathrm{variance}\ \mathrm{of}\ f} dx.
\]
Note that the expectation and variance operators apply to the random
sample $\mathcal{X}_N$. The MISE depends on the bandwidth $\sqrt{t}$
and $f$ in a quite complicated way. The analysis
is simplified when one considers the asymptotic approximation
to the MISE, denoted AMISE, under the consistency requirements that
$t=t_N$ depends on the sample size $N$
such that $t_N\downarrow0$ and $N\sqrt{t_N}\rightarrow\infty$ as
$N\rightarrow\infty$, and $f$ is twice continuously differentiable
\cite{sheather91}.
The asymptotically optimal
bandwidth is then the minimizer of the AMISE. The asymptotic properties
of (\ref{KDE}) under these assumptions are summarized in Appendix \ref
{append:A}.

A key observation about the Gaussian kernel density estimator (\ref
{KDE}) is that
it is the unique solution to the diffusion partial differential
equation (PDE)
%
%
\begin{equation} \label{CH}
\frac{\partial}{\partial t}\hat f(x;t)=\frac
{1}{2}\,\frac{\partial^2}{\partial x^2} \hat f(x;t),\qquad x\in\mathscr
{X}, t>0,
\end{equation}
with $ \mathscr{X}\equiv\mathbb{R} $ and initial condition
$\hat f(x;0)=\Delta(x), $ where
$\Delta(x)=\frac{1}{N}\sum_{i=1}^N\delta(x-X_i)$ is the empirical
density of the data $\mathcal{X}_N$ [here $\delta(x-X_i)$ is the
Dirac measure at~$X_i$].
Equation (\ref{CH}) is the well-known Fourier heat equation
\cite{LarssonThomee}. This link between the Gaussian kernel density
estimator and the Fourier heat equation has been noted in
Chaudhuri and Marron \cite{chaudhuri}. We will, however, go much
further in exploiting this link.
In the heat equation interpretation, the Gaussian kernel in (\ref
{KDE}) is the so-called Green's function \cite{LarssonThomee} for the
diffusion PDE (\ref{CH}).
Thus, the Gaussian kernel density estimator $\hat f(x;t)$ can be
obtained by evolving the solution of the parabolic PDE (\ref{CH}) up
to time $t$.

To illustrate the advantage of the PDE formulation over the more
traditional formulation
(\ref{KDE}), consider the case where
the domain of the data is \textit{known to be} $\mathscr{X}\equiv[0,1]$.
It is difficult to see how (\ref{KDE}) can be easily modified to
account for the finite support of the unknown density. Yet, within the
PDE framework, all we have to do is solve the diffusion equation (\ref
{CH}) over the finite domain $[0,1]$ with
initial condition $\Delta(x)$ and the Neumann boundary condition
\[
\frac{\partial}{\partial x}\hat f(x;t) \Big|_{x=1}=\frac{\partial
}{\partial x}\hat f(x;t) \Big|_{x=0}=0.
\]
The boundary condition ensures that $\frac{d}{dt}\int_\mathscr
{X}\hat f(x;t) \,dx=0$, from where it follows that $\int_\mathscr
{X}\hat f(x;t) \,dx=\int_\mathscr{X}\hat f(x;0) \,dx=1$ for all $t\geq0$.
The analytical solution of this PDE in this case is \cite{bellman}
%
%
\begin{equation}
\label{theta estimator}
\hat f(x;t)=\frac{1}{N}\sum_{i=1}^N\kappa(x,X_i;t),\qquad x\in[0,1],
\end{equation}
%
where
the kernel $\kappa$ is given by
%
%
\begin{equation}
\label{theta func}\quad
\kappa(x,X_i;t) = \sum_{k=-\infty}^\infty
\phi(x,2k+X_i;t)+\phi(x,2k-X_i;t),\qquad x\in[0,1].
\end{equation}
Thus, the kernel accounts for the boundaries in a manner similar to the
boundary correction of the \textit{reflection method} \cite{Silverman}.
We now compare the properties of the
kernel (\ref{theta func}) with the properties of the Gaussian kernel
$\phi$ in (\ref{KDE}).

First, the series representation (\ref{theta func}) is
useful for deriving the small bandwidth properties of the estimator in
(\ref{theta estimator}).
The asymptotic behavior
of $\kappa(x,X_i;t)$ as $t\rightarrow0$ in the interior of the domain
$[0,1]$ is no different from that of the Gaussian kernel, namely,
\[
\sum_{k=-\infty}^\infty
\phi(x,2k+X_i;t)+\phi(x,2k-X_i;t)\sim\phi(x,X_i;t),\qquad t\downarrow0,
\]
for any fixed $x$ in the interior of the domain $[0,1]$. Here
$q(t)\sim z(t), t\downarrow t_0$ stands for $\lim_{t\downarrow
t_0}\frac{q(t)}{z(t)}=1$.
Thus, for small $t$, the estimator (\ref{theta estimator})
behaves like the Gaussian kernel density estimator (\ref{KDE}) in the
interior of $[0,1]$. Near the boundaries
at $x=0,1$, however, the estimator (\ref{theta estimator}) is
consistent, while the Gaussian
kernel density estimator is inconsistent.
In particular, a general result in Appendix \ref{appendix:boundary}
includes as a special case the following boundary property of the
estimator (\ref{theta estimator}):
\[
\mathbb E_f \hat f(x_N;t_N)=f(x_N)+O\bigl(\sqrt{t_N}\bigr),\qquad N\rightarrow\infty,
\]
where $x_N=\alpha t_N$ for some $\alpha\in[0,1]$, and $t_N\downarrow
0$ as $N\rightarrow\infty$. This shows that (\ref{theta estimator})
is consistent at the boundary $x=0$. Similarly, (\ref{theta
estimator}) can be
shown to be consistent at the boundary $x=1$. In contrast,
the Gaussian kernel density estimator (\ref{KDE}) is inconsistent
\cite{WandJones} in the sense that
\[
\mathbb E_f\hat f(0;t_N)=\tfrac{1}{2}f(0)+O\bigl(\sqrt{t_N}\bigr),
\qquad N\rightarrow\infty.
\]
%

The large bandwidth behavior
($t\rightarrow\infty$) of (\ref{theta estimator}) is obtained from
the following equivalent expression for (\ref{theta func}) (see
\cite{bellman}):
%
%
\begin{equation}
\label{theta fourier}
\kappa(x,X_i;t) =\sum_{k=-\infty}^\infty\e^{-k^2\pi^2t/2}\cos
(k\pi x)\cos(k\pi X_i).
\end{equation}
From (\ref{theta fourier}), we immediately see that
%
%
\begin{equation}
\label{theta asympt}
\kappa(x,X_i;t)\sim1+2 \e^{-\pi^2t/2}\cos(\pi x)\cos(\pi X_i) ,\qquad
t\rightarrow\infty, x\in[0,1].
\end{equation}
In other words, as the bandwidth becomes larger and larger, the kernel
(\ref{theta func}) approaches the uniform density on $[0,1]$.
\begin{remark}
An important property of the estimator
(\ref{theta estimator}) is that
the number of local maxima or modes is a nonincreasing function of $t$. This
follows from the \textit{maximum principle} for parabolic
PDE; see, for example, \cite{LarssonThomee}.

For example, a necessary condition for a local maximum at,
say, $(x_0,t_0), t_0>0, x_0\in(0,1)$ is $\frac{\partial^2}{\partial
x^2}\hat f(x_0;t_0)\leq0$. From (\ref{CH}), this implies
$\frac{\partial}{\partial t}\hat f(x_0;t_0)\leq0$, from which it
follows that
there exists an $\varepsilon>0$ such that $\hat f( x_0;t_0)\geq
\hat f(x_0;t_0+\varepsilon)$. As a consequence of this, as $t$ becomes
larger and
larger, the number of local maxima of (\ref{theta estimator}) is a
nonincreasing function. This property is shared by the Gaussian kernel density
estimator (\ref{KDE}) and has been exploited in various ways by
Silverman \cite{Silverman}.
\end{remark}
\begin{example}
Figure \ref{fig:boundary} gives an illustration of the
performance of estimators (\ref{theta estimator})
and (\ref{KDE}), where the true p.d.f. is the beta density
$4(1-x)^3, x\in[0,1]$, and the estimators are build from a
%
%
\begin{figure}

\includegraphics{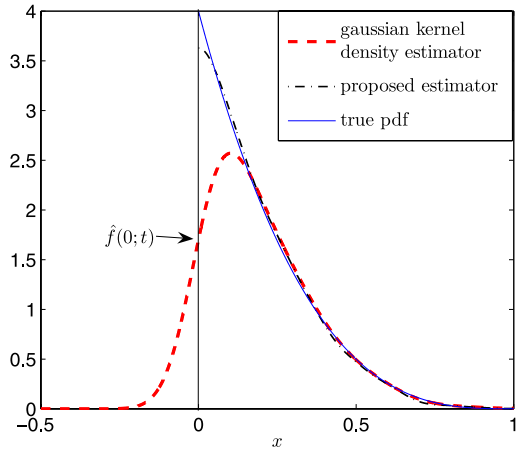}

\caption{Boundary bias in the neighborhood of $x=0$.}
\label{fig:boundary}
\end{figure}
sample of size $N=1000$ with a common bandwidth
$\sqrt{t}=0.05248$.
Note that the Gaussian
kernel density estimator
is close to half the value of the true p.d.f. at the boundary $x=0$.
Overall, the diffusion estimator (\ref{theta estimator}) is much
closer to the true p.d.f. The proposed estimator (\ref{theta estimator})
appears to be the first kernel estimator that does not use boundary
transformation and yet
is consistent at all boundaries
and remains a genuine p.d.f. (is nonnegative
and integrates to one). Existing boundary correction methods \cite
{HallPark02,Alberts,Karunamuni} either
account for the bias at a single end-point, or the resulting estimators
are not genuine p.d.f.'s.
\end{example}
\begin{remark}
In applications such as the smoothed bootstrap
\cite{doufregor01}, there is a need for efficient random variable generation
from the kernel density estimate.
Generation of random variables from the kernel
(\ref{theta func}) is easily
accomplished using the following procedure.
Generate
$ Z\sim\mathsf{N}(0,t)$ and let $Y=X_i+Z$.
Compute $W=Y\operatorname{mod} 2$, and let $X=|W|$. Then it is easy to
show (e.g., using characteristic functions) that $X$ has the density
given by (\ref{theta func}).
\end{remark}

Given the nice boundary bias properties of the estimator that arises as
the solution of the diffusion PDE
(\ref{CH}), it is of interest to investigate if equation (\ref{CH})
can be somehow
modified or generalized to arrive at an even better kernel
estimator. This motivates us to consider in the next section the most
general linear time-homogeneous diffusion PDE as a starting point for
the construction of a better kernel density estimator.

\section{The diffusion estimator}\label{sec3}


Our extension of the simple diffusion model (\ref{CH}) is based on
the smoothing properties of the linear diffusion PDE
%
%
\begin{equation}\label{diffusion} \frac{\partial}{\partial
t}g(x;t)=Lg(x;t),\qquad x\in\mathscr{X}, t>0,
\end{equation}
where the linear differential operator
$L$ is of the form $\frac{1}{2}\frac{d}{d x} (a(x)\frac{d}{d
x} (\frac{\cdot}{p(x)} ) )$, and $a$ and $p$ can be any
arbitrary positive functions on $\mathscr{X}$ with bounded second
derivatives, and the
initial condition is $g(x,0)=\Delta(x)$.
If the set $\mathscr{X}$ is bounded, we add the boundary condition
$\frac{\partial}{\partial x} (\frac{g(x;t)}{p(x)} )=0$
on $\partial\mathscr{X}$, which ensures that the solution
of (\ref{diffusion}) integrates to unity.
The PDE (\ref{diffusion}) describes the p.d.f. of $X_t$ for the It\^o
diffusion process $(X_t,t>0)$ given by \cite{ethier}
%
%
\begin{equation}
\label{Ito}
dX_t= \mu(X_t)\,dt+\sigma(X_t)\,dB_t,
\end{equation}
where the drift coefficient $\mu(x)= \frac{a'(x)}{2 p(x)}$, the
diffusion coefficient $\sigma(x)= \sqrt{\frac{a(x)}{p(x)}}$, the
initial state $X_0$ has distribution $\Delta(x)$, and $(B_t,t>0)$ is
standard Brownian motion.
Obviously, if $a=1$
and $p=1$, we revert to the simpler model (\ref{CH}). What makes the
solution $g(x;t)$ to (\ref{diffusion}) a plausible kernel density
estimator is that $g(x;t)$
is a p.d.f. with the following properties. First, $g(\cdot;0)$ is
identical to the initial condition of (\ref{diffusion}), that is, to
the empirical density $\Delta(x)$. This property is possessed by both
the Gaussian kernel density estimator (\ref{KDE}) and the diffusion
estimator (\ref{theta estimator}).
Second, if $p(x)$ is a p.d.f. on $\mathscr{X}$, then
\[
\lim_{t\rightarrow\infty} g(x;t)=p(x),\qquad x\in\mathscr{X}.
\]
This property is similar to the property that the kernel
(\ref{theta asympt}) and the estimator (\ref{theta estimator})
converge to the uniform density on $\mathscr{X}\equiv[0,1]$ as
$t\rightarrow\infty$. In the context of the diffusion process
governed by (\ref{Ito}), $p$ is the limiting and stationary
density of the diffusion.
Third, similar to the estimator (\ref{theta estimator}) and the
Gaussian kernel density estimator (\ref{KDE}),
we can write the solution of (\ref{diffusion}) as
%
%
\begin{equation}
\label{diff est}
g(x;t)=\frac{1}{N}\sum_{i=1}^N\kappa(x,X_i;t),
\end{equation}
where for each fixed $y\in\mathscr{X}$ the diffusion kernel $\kappa$
satisfies the PDE
%
%
\begin{equation}
\label{FE}
\cases{
\dfrac{\partial}{\partial t} \kappa(x,y;t)= L\kappa(x,y;t),&\quad $x\in
\mathscr{X},
t>0$,\vspace*{2pt}\cr
\kappa(x,y;0) = \delta(x-y), &\quad $x\in\mathscr{X}$.}
\end{equation}

In addition, for each fixed $x\in\mathscr{X}$ the kernel $\kappa$
satisfies the PDE
%
%
\begin{equation}
\label{BE}
\cases{
\dfrac{\partial}{\partial t} \kappa(x,y;t)=L^* \kappa(x,y;t), &\quad
$y\in\mathscr{X},t>0$,\vspace*{2pt}\cr
\kappa(x,y;0) = \delta(x-y), &\quad $y\in\mathscr{X}$,}
\end{equation}
where $L^*$ is of the form
$\frac{1}{2p(y)}\,\frac{\partial}{\partial y} (a(y)\frac{\partial
}{\partial y} (\cdot) )$; that is, $L^*$ is the adjoint operator
of~$L$. Note that $L^*$ is the
\textit{infinitesimal generator} of the It\^o diffusion process
in (\ref{Ito}).
If the set $\mathscr{X}$ has boundaries, we add the Neumann
boundary condition
%
%
\begin{equation}
\label{BC}
\frac{\partial}{\partial x} \biggl(\frac{\kappa(x,y;t)}{p(x)} \biggr)
\bigg|_{x\in\partial\mathscr{X}} =0\qquad \forall t>0
\end{equation}
and $\frac{\partial}{\partial y}\kappa(x,y;t) |_{y\in\partial
\mathscr{X}}=0$ to (\ref{FE}) and (\ref{BE}), respectively. These
boundary conditions ensure that $g(x;t)$ integrates to unity for all
$t\geq0$. The reason that the kernel $\kappa$ satisfies both PDEs
(\ref{FE}) and (\ref{BE})
is that (\ref{FE}) is the Kolmogorov forward equation \cite{ethier}
corresponding to the
diffusion process (\ref{Ito}), and
(\ref{BE}) is a direct consequence of the Kolmogorov backward
equation. We will use the forward and backward equations to derive
the asymptotic properties of the diffusion estimator (\ref{diff est}).
Before we proceed with the asymptotic analysis, we
illustrate how the model (\ref{diffusion}) possesses adaptive
smoothing properties similar to the ones possessed by the adaptive
kernel density estimators \cite{abramson,hall90,hall95,jones94}.
\begin{example}\label{example2}
Suppose that the initial condition of PDE (\ref{diffusion}) is $\Delta
(x)$ with
$N=500\mbox{,}000$ and $X_1,\ldots,X_N$ are independent draws from
$f(x)=1-\cos(6\pi x), x\in[0,1]$. Suppose further that
$p(x)=4(1-x)^3$ and $a(x)=1$ on $[0,1]$. The aim of this example is not
to estimate $f$, but to illustrate the various shapes that the
estimator can take, given data from $f$. Figure
\ref{fig:example_1} shows the solution of the PDE (\ref{diffusion})
for two values of the bandwidth: $\sqrt{t}=4\times10^{-4}$ (small)
and $\sqrt{t}=0.89$ (large).
Since $p(x)$ is the limiting and stationary density of the
%
%
\begin{figure}

\includegraphics{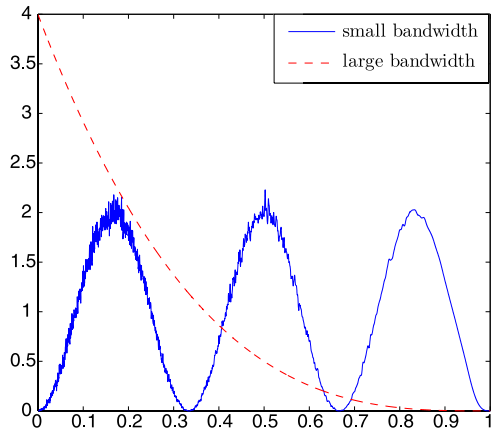}

\caption{Small and large bandwidth behavior of the
diffusion density in Example \protect\ref{example2}.}
\label{fig:example_1}
\end{figure}
diffusion process governed by (\ref{diffusion}), the large
bandwidth density is indistinguishable from $p(x)$. The small
bandwidth density estimate is much closer to $f(x)$ than to $p(x)$. The
crucial feature of the small
bandwidth
density estimate is that $p(x)$ allows for varying degrees
of smoothing across the domain of the data, in particular allowing for
greater smoothing to be
applied in areas of sparse data, and relatively less in the high
density regions. It can be seen from Figure \ref{fig:example_1} that
the small time density estimate is noisier in regions where $p(x)$ is
large (closer to $x=0$), and smoother in regions where $p(x)$ is small
(closer to $x=1$). The adaptive smoothing
is a consequence of the fact that the diffusion kernel (\ref{FE}) has
a state-dependent diffusion coefficient
$\sigma(x)=\sqrt{a(x)/p(x)}$, which helps diffuse the initial density
$\Delta(x)$
at a different rate throughout the state space.
\end{example}
\begin{remark}
Even though there is no analytical expression for the diffusion kernel
satisfying (\ref{FE}), we can write $\kappa$ in terms of a
generalized Fourier series in the case that $\mathscr{X}$ is bounded:
%
%
\begin{equation}
\label{Fourier series}
\kappa(x,y;t)=p(x)\sum_{k=0}^\infty\e^{\lambda_k t}\varphi
_k(x)\varphi_k(y),\qquad x,y\in[0,1],
\end{equation}
where $\{\varphi_k\}$ and $\{\lambda_k\}$ are the eigenfunctions and
eigenvalues of the Sturm--Liouville problem on $[0,1]$:
%
%
\begin{eqnarray}
L^*\varphi_k&=&\lambda_k\varphi_k,\qquad
k=0,1,2,\ldots,\nonumber\\[-8pt]\\[-8pt]
\varphi_k'(0)&=&\varphi_k'(1)=0,\qquad k=0,1,2,\ldots.\nonumber
\end{eqnarray}
It is well known (see, e.g., \cite{LarssonThomee}) that $\{\varphi_k\}
$ forms a complete orthonormal basis with respect to the weight $p$ for
$L^2(0,1)$.
From the expression (\ref{Fourier series}), we can see that
the kernel satisfies the \textit{detailed balance} equation for
a continuous-time Markov process \cite{ethier}
%
%
\begin{equation}
\label{detailed balance}
p(y)\kappa(x,y;t)=p(x)\kappa(y,x;t)\qquad \forall t>0, x,y\in\mathscr{X}.
\end{equation}
The detailed balance equation ensures that the limiting and
stationary density of the diffusion estimator (\ref{diff est}) is
$p(x)$. In addition, the kernel satisfies the
Chapman--Kolmogorov equation
%
%
\begin{equation}
\label{chapman}
\int_\mathscr{X} \kappa(x_1,x_0;t_1) \kappa(x_2,x_1;t_2)
\,dx_1=\kappa(x_2,x_0;t_1+t_2).
\end{equation}

Note that there is no loss of generality in assuming that the domain is
$[0,1]$, because any bounded domain
can be mapped onto $[0,1]$ by a linear transformation.
\end{remark}
%
%
\begin{remark}
When $p(x)$ is a p.d.f.,
an important distance measure between the diffusion estimator (\ref
{diff est}) and $p(x)$ is the divergence measure of Csisz\'ar
\cite{Csiszer}. The Csisz\'ar distance measure between two continuous
probability densities $g$ and $p$ is defined as
\[
\mathcal{D}(g\rightarrow p)=\int_\mathbb{R} p(x) \psi\biggl(\frac
{g(x)}{p(x)} \biggr)\,dx,
\]
where $\psi\dvtx\mathbb{R}_+\rightarrow\mathbb{R}_+$ is a twice
continuously differentiable function;
$\psi(1)=0$; and $\psi''(x)>0$ for all $x\in\mathbb{R}_+$.
The diffusion estimator (\ref{diff est}) possesses the monotonicity property
\[
\frac{d }{d t}\mathcal{D}(g\rightarrow p)=-\frac{1}{2}\int
_{\mathscr{X}} \biggl(\frac{g(x;t)}{p(x)} \biggr)^2\psi'' \biggl(\frac{g(x;t)}{p(x)}
\biggr) \,dx<0,\qquad g\not=p, t>0.
\]
In other words, the distance between the estimator
(\ref{diff est}) and the stationary density $p$ is a monotonically
decreasing function of the bandwidth $\sqrt{t}$. This is why the
solution of (\ref{diffusion}) in Figure
\ref{fig:example_1} approaches $p$ as the bandwidth becomes larger and
larger. Note that Csisz\'ar's family of measures subsumes all of the
information-theoretic distance measures used in practice
\cite{Havrda-Charvat,kapkes87}. For example, if $\psi(x)=\frac{x^\alpha
-x}{\alpha(\alpha-1)}, \alpha\not=0,1$, for some parameter $\alpha
$, then
the family of distances indexed by $\alpha$
includes the \textit{Hellinger distance} for $\alpha=1/2$,
\textit{Pearson's $\chi^2$ discrepancy measure} for
$\alpha=2$, \textit{Neymann's $\chi^2$ measure} for
$\alpha=-1$, the \textit{Kullback--Leibler distance} in the limit as
$\alpha\rightarrow1$ and \textit{Burg's distance} as $\alpha
\rightarrow0$.
\end{remark}

%
%
%
%
%

\section{Bias and variance analysis}\label{sec4}

We now examine the asymptotic bias, variance and MISE of the diffusion
estimator (\ref{diff est}).
In order to derive the asymptotic properties of the proposed estimator,
we need the small bandwidth behavior of
the diffusion kernel satisfying (\ref{FE}).
This is provided by the following lemma.
\begin{lemma}\label{lem:kernel} Assume that the functions $a(x)$ and
$p(x)$ are such that
%
%
\begin{eqnarray}
\label{big assumption}
&\displaystyle c_1=\sqrt{\int_{-\infty}^\infty\biggl(\frac{Lq(z)}{q(z)} \biggr)^2 \,dz}<\infty
,\qquad q(z):=\frac{p(z)}{a^{1/4}(z)p^{1/4}(z)},& \nonumber\\[-8pt]\\[-8pt]
&\displaystyle\lim_{z\rightarrow\infty}\int_{z_0}^z\sqrt{p(s)/a(s)}
\,ds=\infty.&\nonumber
\end{eqnarray}
Then, the leading small bandwidth asymptotic behavior of the kernel
satisfying (\ref{FE}) and (\ref{BE}) on $\mathscr{X}\equiv\mathbb
{R}$ is
\begin{eqnarray*}
\kappa(x,y;t)&\sim&\frac{p(x)}{\sqrt{2\pi t}
[p(x)a(x)a(y)p(y)]^{1/4}}\\
&&{}\times\exp\Biggl\{-\frac{1}{2t} \Biggl[\int_y^x\sqrt{\frac
{p(s)}{a(s)}} \,ds \Biggr]^2 \Biggr\},\qquad t\downarrow0.
\end{eqnarray*}
We denote the asymptotic approximation on the right-hand side by
$\widetilde\kappa(x,y;t)$. Thus, $\kappa(x,y;t)\sim\widetilde
\kappa(x,y;t)$ as $t\downarrow0$.
\end{lemma}

The somewhat lengthy and technical proof is given in Appendix \ref
{appendix:B}. A few remarks about the technical conditions on $a$ and
$p$ now follow. Conditions (\ref{big assumption}) are trivially
satisfied if $a,p$ and its derivatives up to order 2 are all bounded
from above, and $p(x)\geq p_0>0$
and $a(x)\geq a_0>0$. In other words, if we clip $p(x)$ away from zero
and use $a(x)=p^\alpha(x)$ for $\alpha\in[0,1]$, then the conditions
(\ref{big assumption}) are satisfied. Such clipping procedures have
been applied in the traditional kernel density estimation setting, see
\cite{abramson,choi,hall95,hall02,jones94}.
Note that the conditions are more easily satisfied when $p$ is
heavy-tailed. For example,
if $a(x)=p(x)$, then $p$ could be any regularly varying p.d.f. of the form
$p\varpropto(1+|x|)^{-\alpha}, \alpha>1 $. Lemma \ref{lem:kernel}
is required for deriving the asymptotic properties of the estimator,
all collected in the following theorem.
\begin{theorem}
\label{diffusion theorem}
Let $t=t_N$ be such that $\lim_{N\rightarrow\infty}t_N=0$, $\lim
_{N\rightarrow\infty} N\sqrt{t_N}=\infty$. Assume that $f$ is twice
continuously differentiable and that the domain $\mathscr{X}\equiv
\mathbb{R}$. Then:
\begin{enumerate}
\item
The pointwise bias has the asymptotic behavior
%
%
\begin{equation}
\label{theorem:bias}
\mathbb E_f [g(x;t)]-f(x)= t Lf(x)+O(t^2),\qquad N\rightarrow\infty.
\end{equation}

\item
The integrated squared bias has the asymptotic behavior
%
%
\begin{equation} \label{diffusion_bias}
\Vert\mathbb E_f[g(\cdot
;t)]-f\Vert^2\sim t^2 \Vert Lf\Vert^2=\tfrac{1}{4}t^2\bigl\Vert\bigl(a
(f/p)'\bigr)'\bigr\Vert^2,\qquad
N\rightarrow\infty.
\end{equation}

\item The pointwise variance has the asymptotic behavior
%
%
\begin{equation}
\label{theorem:var}
\operatorname{Var}_f[g(x;t)]\sim\frac{f(x)}{2N\sqrt{\pi t}
\sigma(x)},\qquad
N\rightarrow\infty,
\end{equation}
where $\sigma^2(x)=a(x)/p(x)$.
\item The integrated variance has the asymptotic behavior
%
%
\begin{equation} \label{diffusion_var}
\int\operatorname{Var}_f[g(x;t)] \,dx\sim\frac{\mathbb E_f[\sigma
^{-1}(X)]}{2N\sqrt{\pi t}},\qquad N\rightarrow\infty.
\end{equation}
%
%
\item Combining the leading order bias and variance terms gives the
asymptotic approximation to the MISE
%
%
\begin{equation}
\operatorname{AMISE}\{g\}(t)=\frac{1}{4}t^2\bigl\Vert\bigl(a(f/p)'\bigr)'\bigr\Vert^2+\frac
{\mathbb E
_f[\sigma^{-1}(X)]}{2N\sqrt{\pi t}}.
\end{equation}
\item Hence, the square of the asymptotically optimal bandwidth is
%
%
\begin{equation}\label{diffusion_t*}
t^*= \biggl(\frac{\mathbb E_f[\sigma^{-1}(X)]}{2N\sqrt{\pi} \Vert Lf \Vert^2} \biggr)^{2/5},
\end{equation}
which gives the minimum
%
%
\begin{equation}\label{diffusion_AMISE}
\min_t\operatorname{AMISE}\{g\}(t)=N^{-4/5} \frac{5
[\mathbb E
_f\sigma^{-1}(X)]^{4/5} \Vert Lf \Vert^{2/5}}{2^{14/5}\pi^{2/5}} .
\end{equation}
\end{enumerate}
\end{theorem}

The proof is given in Appendix \ref{appendix:C}.

We make the following observations.
First, if $p\not\equiv f$, the rate of convergence of
(\ref{diffusion_AMISE}) is $O(N^{-4/5})$, the same as the rate of the
Gaussian kernel density estimator in (\ref{rate}).
The multiplicative constant of $N^{-4/5}$ in (\ref{diffusion_AMISE}),
however, can be made very small by choosing $p$ to be a
\textit{pilot density estimate} of $f$.
Preliminary or pilot density estimates are used in most adaptive kernel
methods \cite{WandJones}.
Second, if $p\equiv
f$, then the leading bias term (\ref{theorem:bias}) is 0.
In fact, if $f$ is infinitely smooth, the pointwise bias
is exactly zero, as can be seen from
\[
\mathbb E_f [g(x;t)]= \sum_{k=0}^\infty\frac{t^k}{k!} L^kf(x) ,\qquad f\in
C^\infty,
\]
where $L^{n+1}=LL^n$ and $L^0$ is the identity operator.
In addition, if $a=p\varpropto1$, then the bias term
(\ref{theorem:bias}) is equivalent to the bias term (\ref{kde_bias})
of the Gaussian kernel density estimator.
Third, (\ref{theorem:var}) suggests that in regions where the pilot density
$p(x)$ is large [which is equivalent to small diffusion coefficient
$\sigma(x)$] and
$f(x)$ is large, the pointwise variance will be large. Conversely, in
regions with few observations [i.e., where the diffusion coefficient
$\sigma(x)$
is high and $f(x)$ is small] the pointwise variance is low. In other words,
the ideal variance behavior results when the diffusivity $\sigma(x)$
behaves inversely proportional to $f(x)$.

\subsection{Special cases of the diffusion estimator}\label{sec41}

We shall now show that the diffusion kernel estimator (\ref{diff est})
is a generalization of some well-known modifications of the Gaussian
kernel density estimator (\ref{KDE}). Examples of modifications and
improvements subsumed as special cases of (\ref{diff est}) are as
follows.

\begin{enumerate}
\item If $a(x)=p(x)\varpropto1$ in (\ref{diff est}) and $\mathscr
{X}\equiv\mathbb{R}$, then the kernel $\kappa$ reduces to the
Gaussian kernel and we obtain (\ref{KDE}).
\item If $a(x)=1$ and $p(x)= f_p(x)$, where $f_p$ is a clipped pilot
density estimate of $f$ (see \cite{abramson,hall02,jones94}), then
from Lemma \ref{lem:kernel}, we have
\[
\kappa(x,y;t)\sim\widetilde\kappa(x,y;t)=\frac{f_p(x)}{\sqrt{2\pi
t} (f_p(x)f_p(y))^{1/4}}\exp\biggl\{-\frac{1}{2t} \biggl[\int_y^x\sqrt{f_p(s)}
\,ds \biggr]^2 \biggr\}.
\]
Thus, in the neighborhood of $y$ such that $|x-y|=O(t^\beta), \beta
>1/3$, we have
\[
\kappa(x,y;t)\sim\frac{1}{\sqrt{2\pi t/f_p(x)}}\exp\biggl\{-\frac
{(x-y)^2}{2t/f_p(x)} \biggr\},\qquad t\downarrow0.
\]
In other words, in the neighborhood of $y$, $\kappa$ is
asymptotically equivalent to a Gaussian kernel with mean
$y$ and bandwidth $\sqrt{t/f_p(y)}$, which is precisely
the Abramson's variable bandwidth \cite{abramson}
modification as applied to the Gaussian kernel.
Abramson's square root law states
that the asymptotically optimal variable bandwidth
is proportional to $ f_p^{-1/2}(y)$.

%
\item If we choose $a(x)=p(x)= f_p(x)$, then in an $O(t^\beta), \beta
>0$ neighborhood of $y$, the kernel $\kappa(x,y;t)$ behaves
asymptotically as a Gaussian kernel with location $y+\frac{t}{2}\frac
{ f_p'(y)}{f_p(y)}$ and bandwidth $\sqrt{t}$:
\[
\kappa(x,y;t)\sim\frac{1}{\sqrt{2\pi t}}\exp\biggl\{-\frac{1}{2t}
\biggl(x-y-\frac{t}{2}\frac{f_p'(y)}{f_p(y)} \biggr)^2 \biggr\} ,\qquad
t\downarrow0.
\]
This is precisely the data sharpening
modification described in \cite{Samiuddin}, where
the locations of the data points\vspace*{1pt} are shifted prior to the
application of the kernel density estimate. Thus, in our paradigm,\vspace*{1pt} data
sharpening is equivalent to using the diffusion (\ref{diffusion}) with drift
$\mu(x)=\frac{f_p'(x)}{2f_p(x)}$ and diffusion coefficient
$\sigma(x)=1$.

\item Finally, if we set $p(x)= f_p(x)$ and $a(x)=p^{\alpha}(x)$,
$\alpha\in[0,1]$, then we obtain a method that is a combination of
both the data sharpening and the variable bandwidth of Abramson. The
kernel $\kappa$
behaves asymptotically [in an $O(t^\beta), \beta>1/3$ neighborhood of
$y$] like a Gaussian kernel with location $y+t\mu(y)=y+\frac{\alpha
t}{2} f_p^{\alpha-2}(y)f_p'(y)$ and bandwidth $\sqrt{t \sigma
^2(y)}=\sqrt{t f_p^{\alpha-1}(y)}$. Similar variable location and
scale kernel density estimators are considered in \cite{jones94}.

\end{enumerate}
The proposed method thus unifies many of the already existing ideas for
variable scale and location kernel density estimators. Note that these
estimators all have one common feature: they compute a pilot density
estimate (which is an infinite-dimensional parameter) prior to the main
estimation step.

Our choice for $a(x)$ will be motivated by regularity properties of the
diffusion process
underlying the smoothing kernel. In short, we prefer to choose $a(x)=1$
so as to
make the diffusion process in (\ref{Ito}) nonexplosive with a
well-defined limiting distribution.
A necessary and sufficient condition for explosions is
Feller's test~\cite{feller}.
\begin{theorem}[(Feller's test)]
Let $\mu(x)>0$ and $\sigma(x)>0$ be bounded and continuous.
Then the diffusion process (\ref{Ito}) explodes if and only if there
exists $z\in\mathbb{R}$ such that either one of the following
two conditions holds:
\begin{enumerate}
\item
\[
\int_{-\infty}^{z} \int_x^z \exp\biggl( \int_{x}^y \frac{2\mu
(s)}{\sigma^2(s)} \,d s \biggr)\sigma^{-2}(y)\,d y \,d x<\infty,
\]
\item
\[
\int_{z}^{\infty} \int_z^x \exp\biggl( \int_{x}^y \frac{2\mu
(s)}{\sigma^2(s)} \,d s \biggr)\sigma^{-2}(y)\,d y\, d x<\infty.
\]
\end{enumerate}
\end{theorem}

A corollary of Feller's test is that when $\mu(x)=0$ both of Feller's
conditions fail, and diffusions of the form
$d X_t=\sigma_t \,d W_t$ are nonexplosive.

Since in our case we have
$\sigma^2(x)=a(x)/p(x)$ and $a(x)=\exp( \int_{x_0}^x 2\mu(y)/\break\sigma
^2(y) \,d y )$, Feller's condition becomes the following.
\begin{proposition}[(Feller's test)]\label{prop1}
Given $a(x)$ and $p(x)$
in (\ref{diffusion}), the diffusion process (\ref{Ito})
explodes if and only if there exists $z\in\mathbb{R}$ such that
either one of the following
two conditions holds:
\begin{enumerate}
\item
\[
\int_{-\infty}^{z} \int_x^z \frac{p(y)}{a(x)}\,d y \,d x<\infty,
\]
\item
\[
\int_{z}^{\infty} \int_z^x \frac{p(y)}{a(x)}\,d y \,d x<\infty.
\]
\end{enumerate}
\end{proposition}

The easiest way to ensure nonexplosiveness of the underlying
diffusion process and the existence of a limiting distribution is to
set $a(x)=1$, which corresponds to $\mu(x)=0$. Note that a necessary
condition for the
existence of a limiting p.d.f. is the existence of $z$ such that $\int
_z^\infty1/a(x)\,d x=\infty$.
In this case, both of Feller's conditions fail. The nonexplosiveness
property ensures that generation of random variables from the
diffusion estimator does not pose any technical problems.

\section{Bandwidth selection algorithm}\label{sec5}

Before we explain how to estimate the bandwidth $\sqrt{t^*}$ in (\ref
{diffusion_t*}) of the diffusion estimator (\ref{diff est}), we explain
how to estimate the bandwidth $\sqrt{_*t}$ in (\ref{Gaussian kde
bandwidth}) (see Appendix \ref{append:A}) of the Gaussian kernel
density estimator (\ref{KDE}). Here, we present a new plug-in
bandwidth selection procedure
based on the ideas in \cite{sheather96,JonesMarronPark,Marron,sheather91}
to achieve unparalleled practical performance. The
highlighting feature of the proposed method is that it
does not use normal reference rules and is thus
completely data-driven.

It is clear from (\ref{Gaussian kde bandwidth}) in Appendix \ref
{append:A} that to compute the optimal $_*t$ for the Gaussian
kernel density estimator (\ref{KDE}) one needs to estimate
the functional $\Vert f''\Vert^2$. Thus, we consider the problem of estimating
$\Vert f^{(j)}\Vert^2$ for an arbitrary integer
$j\geq1$.
The identity $\Vert f^{(j)}\Vert^2=(-1)^j\mathbb E_f[f^{(2j)}(X)]$
suggests two
possible plug-in estimators. The first one is
%
%
\begin{eqnarray}
\label{gaussian plug-in estimator}
(-1)^j\widehat{\mathbb E_ff^{(2j)}}:\!&=&\frac{(-1)^j}{N}\sum_{k=1}^N
\hat
f^{(2j)}(X_k;t_j)\nonumber\\[-8pt]\\[-8pt]
&=&\frac{(-1)^j}{N^2}\sum_{k=1}^N\sum_{m=1}^N \phi
^{(2j)}(X_k,X_m;t_j),\nonumber
\end{eqnarray}
where $\hat f$ is the Gaussian kernel density estimator (\ref{KDE}).
The second estimator is
%
%
\begin{eqnarray}
\label{gaussian plug-in estimator 2}
\bigl\Vert \widehat{f^{(j)}}\bigr\Vert^2:\!&=&\bigl\Vert\hat f^{(j)}(\cdot;t)\bigr\Vert^2\nonumber\\
&=&\frac
{1}{N^2}\sum_{k=1}^N\sum_{m=1}^N \int_\mathbb{R}\phi
^{(j)}(x,X_k;t_j)\phi^{(j)}(x,X_m;t_j) \,dx\\
&=&\frac{(-1)^j}{N^2}\sum_{k=1}^N\sum_{m=1}^N
\phi^{(2j)}(X_k,X_m;2 t_j),\nonumber
\end{eqnarray}
where
the last line is a simplification following easily from the fact that the
Gaussian kernel $\phi$ satisfies the Chapman--Kolmogorov
equation (\ref{chapman}). For a given bandwidth, both estimators
$(-1)^j\widehat{\mathbb E_ff^{(2j)}}$ and
$\widehat{\Vert f^{(j)}\Vert}{}^2$ aim to estimate the same quantity, namely
$\Vert f^{(j)}\Vert^2$. We select $t_j$ so that both estimators (\ref
{gaussian plug-in estimator}) and
(\ref{gaussian plug-in estimator 2}) are asymptotically
equivalent in the mean square error sense. In other words,
we choose $t_j={_*}t_j$ so that both
$(-1)^j\widehat{\mathbb E_ff^{(2j)}}$ and
$\widehat{\Vert f^{(j)}\Vert}{}^2$ have equal asymptotic mean
square error. This gives the following proposition.
\begin{proposition}\label{prop2}
The estimators
$(-1)^j\widehat{\mathbb E_ff^{(2j)}}$ and $\widehat{\Vert f^{(j)}\Vert}{}^2$
have the same asymptotic mean square error when
%
%
\begin{equation}
\label{l-stage}
_*t_j= \biggl( \frac{1+{1}/{2^{j+1/2}}}{3} \frac{1\times3\times
5\times\cdots\times(2j-1)}{N\sqrt{\pi/2}
\Vert f^{(j+1)}\Vert^2} \biggr)^{{2}/({3+2j})}.
\end{equation}
\end{proposition}
\begin{pf}
The arguments are similar to
the ones used in \cite{WandJones}. Under the assumptions that $t_j$
depends on $N$ such that\vspace*{1pt}
$\lim_{N\rightarrow\infty}t_j=0$ and\break
$\lim_{N\rightarrow\infty}N t_j^{j+1/2}=\infty$, we can take the
expectation of the estimator
(\ref{gaussian plug-in estimator}) and obtain the expansion $(t_j=t)$:
\begin{eqnarray*}
&&\mathbb E_f \bigl[\widehat{\mathbb E_ff^{(2j)}} \bigr] \\
&&\qquad =\frac{1}{N}\phi^{(2j)}(0,0;t)+\frac{N-1}{N}\iint f(x)f(y)\phi
^{(2j)}(x,y;t) \,dx\,dy\\
&&\qquad =-\frac{1\times3\times\cdots\times(2j-1)}{ t^{j+1/2} \sqrt{2\pi
} N}\\
&&\qquad\quad{}+ \int f(x) \biggl(f^{(2j)}(x)+ \frac{t}{2}f^{2(j+1)}(x)+o(t) \biggr)\,dx
+O(N^{-1})\\
&&\qquad=-\frac{1\times3\times5\times\cdots\times(2j-1)}{ t^{j+1/2}
\sqrt{2\pi} N}+ \frac{t}{2}\bigl\Vert f^{(j+1)}\bigr\Vert^2\\
&&\qquad\quad{} +(-1)^j\bigl\Vert f^{(j)}\bigr\Vert^{2}+O(N^{-1}),\qquad N\rightarrow\infty.
\end{eqnarray*}
Hence, the squared bias has asymptotic behavior ($N\rightarrow\infty$)
\[
\bigl((-1)^j\mathbb E_f \bigl[\widehat{\mathbb E_ff^{(2j)}} \bigr]-\bigl\Vert f^{(j)}\bigr\Vert^2
\bigr)^2\sim
\biggl(\frac{1\times3\times\cdots\times(2j-1)}{ t^{j+1/2} \sqrt{2\pi}
N}- \frac{t}{2}\bigl\Vert f^{(j+1)}\bigr\Vert^2 \biggr)^2.
\]
A similar argument (see \cite{WandJones}) shows that
the variance is of the order $O(N^{-2}\times\break t^{-2j-1/2})$, which is of lesser
order than the squared bias. This implies that the leading order term
in the asymptotic mean square error of $\widehat{\mathbb E_ff^{(2j)}}$ is
given by the asymptotic squared bias.
There is no\vspace*{1pt} need to derive the asymptotic expansion of
$\mathbb E_f [\widehat{\Vert f^{(j)}\Vert}{}^2 ]$, because
inspection of (\ref{gaussian plug-in estimator 2})
and (\ref{gaussian plug-in estimator}) shows that
$\widehat{\Vert f^{(j)}\Vert}{}^2$ exactly equals $(-1)^j\widehat{\mathbb E
_ff^{(2j)}}$ when the latter is evaluated at $2 t_j$. In other words,
\begin{eqnarray*}
(-1)^j\mathbb E_f \bigl[\widehat{\bigl\Vert f^{(j)}\bigr\Vert}^2 \bigr]&=&-\frac{1\times
3\times
5\times\cdots\times(2j-1)}{ (2t)^{j+1/2} \sqrt{2\pi} N}\\
&&{}+
t\bigl\Vert f^{(j+1)}\bigr\Vert^2+ O(1+N^{-1}).
\end{eqnarray*}
Again, the leading term of the asymptotic mean square error of
$\widehat{\Vert f^{(j)}\Vert}{}^2$ is given by the leading term of the squared
bias of $\widehat{\Vert f^{(j)}\Vert}{}^2$.
Thus, equalizing the asymptotic mean squared error of both estimators
is the same as equalizing their respective asymptotic squared biases.
This yields the equation
\begin{eqnarray*}
&&\biggl(\frac{1\times3\times\cdots\times(2j-1)}{ (2t)^{j+1/2} \sqrt{2\pi
} N}- t\bigl\Vert f^{(j+1)}\bigr\Vert^2 \biggr)^2 \\
&&\qquad= \biggl(\frac{1\times3\times\cdots\times
(2j-1)}{ t^{j+1/2} \sqrt{2\pi} N}- \frac{t}{2}\bigl\Vert f^{(j+1)}\bigr\Vert^2 \biggr)^2.
\end{eqnarray*}
The positive solution of the equation yields the desired ${_*}t_j$.
\end{pf}

Thus, for example,
%
%
\begin{equation}
\label{_*t_2}
_*t_2= \biggl(\frac{8+\sqrt{2}}{24} \frac{3}{N\sqrt{\pi/2}
\Vert f^{(3)}\Vert^2} \biggr)^{2/7}
\end{equation}
is our bandwidth choice for the estimation of $\Vert f''\Vert^2$.
We estimate each $_*t_j$ by
%
%
\begin{equation} \label{l-stage estimate}
{_*}\hat{t}_j= \biggl(\frac{1+{1}/{2^{j+1/2}}}{3} \frac{
1\times3\times5\times\cdots\times(2j-1)}{N\sqrt{\pi/2}
\widehat{\Vert f^{(j+1)}\Vert}{}^2} \biggr)^{{2}/({3+2j})}.
\end{equation}
%
Computation of
$\widehat{\Vert f^{(j+1)}\Vert}{}^2$ requires
estimation of ${_*}t_{j+1}$ itself, which in turn requires
estimation of ${_*}t_{j+2}$, and so on, as seen from formulas (\ref
{gaussian plug-in estimator 2}) and (\ref{l-stage estimate}). We are
faced with the problem of estimating the infinite
sequence $ \{{_*}t_{j+k},k\geq1 \}$. It is clear, however, that given
$_*t_{l+1}$ for some $l>0$ we can estimate all $\{ _*t_j, 1\leq j \leq
l\}$ recursively, and then estimate $_* t$ itself from
(\ref{Gaussian kde bandwidth}).
This motivates the \textit{$l$-stage direct plug-in bandwidth selector}
\cite{sheather96,sheather91,WandJones},
defined as follows.
\begin{enumerate}
\item For a given integer $l>0$, estimate $_*t_{l+1}$ via (\ref
{l-stage}) and $\Vert f^{(l+2)}\Vert^2$ computed by assuming that $f$
is a
normal density with mean and variance estimated from the data. Denote
the estimate by $_*\hat{t}_{l+1}$.\vspace*{1pt}
\item Use $_*\hat{t}_{l+1}$ to estimate $\Vert f^{(l+1)}\Vert^2$ via the
plug-in estimator (\ref{gaussian plug-in estimator 2}) and $_*\hat
{t}_{l}$ via
(\ref{l-stage estimate}). Then use $_*\hat{t}_{l}$ to
estimate $_*\hat{t}_{l-1}$ and so on until we obtain an estimate of
$_*\hat{t}_{2}$.
\item Use the estimate of $_*\hat{t}_{2}$ to compute
$_*\hat{t}$ from (\ref{Gaussian kde bandwidth}).
\end{enumerate}
The $l$-stage direct plug-in bandwidth selector thus involves the
estimation of $l$ functionals $\{\Vert f^{(j)}\Vert,2\leq j\leq l+1\}$
via the plug-in
estimator (\ref{gaussian plug-in estimator 2}).
We can describe the procedure in a more abstract
way as follows. Denote the functional dependence of $_*\hat t_j$ on
$_*\hat t_{j+1}$ in formula (\ref{l-stage estimate}) as
\[
_*\hat t_j=\gamma_j(_*\hat t_{j+1}).
\]
It is then clear that $_*\hat t_j=\gamma_j(\gamma_{j+1}(_*\hat
t_{j+2}))=\gamma_j(\gamma_{j+1}(\gamma_{j+2}(_*\hat t_{j+3})))=\cdots.$
For simplicity of notation, we define the composition
\[
\gamma^{[k]}(t)=\underbrace{\gamma_1( \cdots\gamma_{k-1}(\gamma_k}_{k
\ \mathrm{times}}(t)) \cdots),\qquad k\geq1.
\]
%
Inspection of formulas (\ref{l-stage estimate}) and (\ref{Gaussian
kde bandwidth})
shows that the estimate of $_* t$ satisfies
\begin{eqnarray*}
_*\hat t&=&\xi{_*}\hat t_1=\xi\gamma^{[1]}(_*\hat t_{2})=\xi\gamma
^{[2]}({_*}\hat t_{3})=\cdots=\xi\gamma^{[l]}({_*}\hat t_{1+l}),\\
\xi &=& \biggl( \frac{6\sqrt{2}-3}{7} \biggr)^{2/5}\approx0.90.
\end{eqnarray*}
Then, for a given integer $l>0$, the $l$-stage direct plug-in bandwidth
selector consists of computing
\[
_*\hat t=\xi\gamma^{[l]}(_* t_{l+1}),
\]
where $_*t_{l+1}$ is estimated via (\ref{l-stage}) by assuming that
$f$ in $\Vert f^{(l+2)}\Vert^2$ is a normal density with mean and variance
estimated from the data.
The weakest point of this procedure is that we assume that the true $f$
is a Gaussian density in order to compute
$\Vert f^{(l+2)}\Vert^2$. This assumption can lead to arbitrarily bad
estimates of $_*{t}$, when, for example, the true $f$ is far from being
Gaussian. Instead, we propose to find a solution to the nonlinear equation
%
%
\begin{equation}
t=\xi\gamma^{[l]}(t),
\end{equation}
for some $l$,
using either fixed point iteration or Newton's method with initial
guess $t=0$. The fixed point iteration version is formalized in the
following algorithm.
%
\begin{algorithm}[(Improved Sheather--Jones)]
Given $l>2$, execute the following steps:
\label{my plug-in}
\begin{enumerate}
\item initialize with $z_0=\varepsilon$, where $\varepsilon$ is machine
precision, and $n=0$;
\item set $z_{n+1}=\xi\gamma^{[l]}(z_n)$;
\item if $|z_{n+1}-z_n|<\varepsilon$, stop and set $_*\hat t=z_{n+1}$;
otherwise, set $n:=n+1$ and repeat from step 2;
\item deliver the Gaussian kernel density estimator (\ref{KDE})
evaluated at $_*\hat t$ as the final estimator of $f$, and $_*\hat
t_2=\gamma^{[l-1]}(z_{n+1})$ as the bandwidth
for the optimal estimation of $\Vert f''\Vert^2$.
\end{enumerate}
\end{algorithm}

Numerical experience suggests the following. First, the
fixed-point algorithm does not fail to find a root of
the equation $t=\xi\gamma^{[l]}(t)$. Second, the
root appears to be unique. Third, the solutions to the equations
\[
t=\xi\gamma^{[5]}(t)
\]
and
\[
t=\xi\gamma^{[l+5]}(t)
\]
for any $l>0$
do not differ in any practically meaningful way. In other words, there
were no gains to be had by increasing the stages of the bandwidth
selection rule beyond $l=5$. We recommend setting $l=5$.
Finally, the numerical
procedure for the computation of $\gamma^{[5]}(t)$ is fast when
implemented using the
Discrete Cosine Transform \cite{Matlab}. 

The plug-in method described in Algorithm \ref{my plug-in} has
superior practical performance compared to existing
plug-in implementations, including the particular
\textit{solve-the-equation} rule of Sheather and Jones
\cite{sheather91,WandJones}. Since we borrow many of the fruitful ideas
described in \cite{sheather91} (which in turn build upon the work of
Hall, Park and Marron \cite{hall87,park1990}), we call our new
algorithm the Improved Sheather--Jones (ISJ) method.

To illustrate the significant improvement of the plug-in method in
Algorithm \ref{my plug-in}, consider, for example,
the case where $f$ is a mixture of two Gaussian densities with a common
variance of $1$ and means of $-30$ and $30$.

%
\begin{figure}

\includegraphics{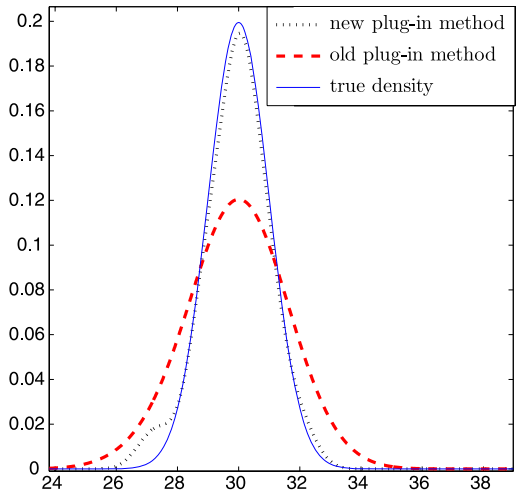}

\caption{The Improved Sheather--Jones bandwidth selection
rule in Algorithm \protect\ref{my plug-in} leads to improved performance
compared to the original plug-in rule that uses the normal reference
rule.} \label{fig:SJ}
\end{figure}

Figure \ref{fig:SJ} shows the right mode of $f$, and the
two estimates resulting from the old plug-in rule
\cite{sheather91} and the plug-in rule of Algorithm
\ref{my plug-in}. The left mode is not displayed, but looks similar.
The integrated squared error using the
new plug-in bandwidth estimate, $\Vert f-\hat f(\cdot;{_*}\hat t)\Vert
^2$, is
one 10th of the error using the old
bandwidth selection rule.

\subsection{Experiments with normal reference rules}\label{sec51}

The result of Figure \ref{fig:SJ}
is not an isolated case, in which the normal reference rules do not
perform well. We performed a comprehensive simulation study in order to
compare the Improved Sheather--Jones (ISJ) (Algorithm \ref{my plug-in})
with the original
(vanilla) Sheather--Jones (SJ) algorithm \cite{sheather91,WandJones}.

%
\begin{table}
\tabcolsep=0pt
\caption{Results over 10 independent simulation experiments. In all
cases the domain was assumed to be $\mathbb{R}$. Many
test problems are taken from \protect\cite{MarronWand}. In the table $\mathsf
{N}(\mu,\sigma^2)$, denotes a Gaussian density with mean $\mu$ and
variance $\sigma^2$}\label{ISJ table}
\begin{tabular*}{\tablewidth}{@{\extracolsep{\fill}}lccc@{}}
\hline
\textbf{Case} & \textbf{Target density} $\bolds{f(x)}$ & $\bolds{N}$ & \textbf{Ratio}\\
\hline
1 (claw)&$\frac{1}{2}\mathsf{N}(0,1)+\sum_{k=0}^4\frac
{1}{10}\mathsf{N} (\frac{k}{2}-1, (\frac{1}{10} )^{2} )$ &
$10^3$ & 0.72\\
& & $10^4$ & 0.94\\
[4pt]
2 (strongly skewed) & $\sum_{k=0}^7\frac{1}{8}\mathsf{N} (3 ( (\frac
{2}{3} )^k-1 ), (\frac{2}{3} )^{2k} )$ &
$10^3$ & 0.69\\
& & $10^4$ & 0.84\\
[4pt]
3 (kurtotic unimodal)&$\frac{2}{3}\mathsf{N}(0,1)+\frac{1}{3}\mathsf
{N} (0, (\frac{1}{10} )^2 )$
& $10^2$ & 0.78\\
& & $10^3$ & 0.93\\
[4pt]
4 (double claw)& \multicolumn{1}{c}{$\frac{49}{100}\mathsf{N} (-1, (\frac
{2}{3} )^2 )+\frac{49}{100}\mathsf{N} (1, (\frac{2}{3} )^2 )$} & $10^5$ & 0.35
\\[4pt]
& \multicolumn{1}{c}{\hspace*{-16.4pt}\hspace*{48pt}$+\,\frac
{1}{350}\sum_{k=0}^6\mathsf{N} (\frac{k-3}{2}, (\frac{1}{100} )^2
)$} & $10^6$ & 0.10\\
[4pt]
5 (discrete comb)& $\frac{2}{7}\sum_{k=0}^2\mathsf
{N} (\frac{12k-15}{7}, (\frac{2}{7} )^2 )+\frac{1}{21}\sum
_{k=8}^{10} \mathsf{N} (\frac{2k}{7}, (\frac{1}{21} )^2 ) $
& $10^3$ & 0.45\\
& & $10^4$ & 0.27\\
[4pt]
6 (asymmetric &\multicolumn{1}{c}{$\frac{46}{100}\sum_{k=0}^1
\mathsf{N} (2k-1, (\frac{2}{3} )^2 )+
\sum_{k=1}^3\frac{1}{300}\mathsf{N} (-\frac{k}{2}, (\frac{1}{100}
)^2 )$}  & $10^4$ & 0.68\\[4pt]
\quad double claw) & \multicolumn{1}{c}{\hspace*{-52.43pt}\hspace*{48pt}$+\,\sum_{k=1}^3\frac{7}{300}\mathsf{N} (\frac{k}{2}, (\frac
{7}{100} )^2 )$\hspace*{52.43pt}} & $10^6$ & 0.24\\
[4pt]
7 (outlier)&$ \frac{1}{10}\mathsf{N}(0,1)+\frac{9}{10}\mathsf{N}
(0, (\frac{1}{10} )^2 ) $ & $10^3$ & 1.01\\
& & $10^5$ & 1.00\\
[4pt]
8 (separated bimodal)& $\frac{1}{2}\mathsf{N}(-12,\frac{1}{4})+\frac
{1}{2}\mathsf{N}(12,\frac{1}{4})$ & $10^2$ & 0.33\\
& & $10^3$ & 0.64\\
[4pt]
9 (skewed bimodal)&$\frac{3}{4}\mathsf{N}(0,1)+\frac{1}{4}\mathsf
{N} (\frac{3}{2}, (\frac{1}{3} )^2 )$ & $10^3$ & 1.02\\
& & $10^4$ & 1.00\\
[4pt]
10 (bimodal)&$\frac{1}{2}\mathsf{N} (0, (\frac{1}{10} )^2 )+\frac
{1}{2}\mathsf{N} (5,1 )$ & $10^2$ & 0.31\\
& & $10^3$ & 0.70\\
[4pt]
11 &Log-Normal with $\mu=0$ and $\sigma=1$ & $10^3$ & 0.82\\
& & $10^4$ & 0.80\\
[4pt]
12 (asymmetric claw)&$ \frac{1}{2}\mathsf{N}(0,1)+\sum_{k=-2}^2\frac
{2^{1-k}}{31}\mathsf{N} (k+\frac{1}{2}, (\frac{2^{-k}}{10} )^2 ) $ &
$10^3$ & 0.76\\
& & $10^4$ & 0.59 \\
[4pt]
13 (trimodal)&$ \frac{1}{3}\sum_{k=0}^2\mathsf{N}(80k;(k+1)^4) $ &
$10^2$ & 0.21\\
& & $10^3$ & 0.17\\
[4pt]
14 (5-modes)& $\frac{1}{5}\sum_{k=0}^4\mathsf{N}(80k;(k+1)^2)$ &
$10^3$ & 0.07\\
& & $10^4$ & 0.18\\
[4pt]
15 (10-modes)&
$\frac{1}{10}\sum_{k=0}^9\mathsf{N}(100k;(k+1)^2)$
& $10^3$ & 0.12\\
& & $10^4$ & 0.07\\
[4pt]
16 (smooth comb)& $\sum_{k=0}^5 \frac{2^{5-k}}{63}\mathsf{N}(
\frac{65-96/2^k}{21};\frac{(32/63)^2}{2^{2k}} )$ &
$10^4$ & 0.40\\
& & $10^5$ & 0.34\\
\hline
\end{tabular*}
\end{table}

Table \ref{ISJ table} shows
the average results over 10 independent trials for a number of different
test cases. The second column displays the target density and the
third column shows the sample size used for the experiments. The last column
shows our criterion for comparison:
\[
R=\frac{\Vert\hat f(\cdot;{_*}\hat t)-f\Vert^2 }{\Vert\hat f(\cdot;
t_{\mathrm{SJ}})-f\Vert^2},
\]
that is, the ratio of the integrated squared error of the new ISJ
estimator to the integrated squared error of the original SJ estimator.
Here, $t_{\mathrm{SJ}}$ is the bandwidth computed
using the original Sheather--Jones method \cite{sheather91,WandJones}.

The results in Table \ref{ISJ table} show that the improvement in the
integrated squared error
can be as much as ten-fold, and the ISJ method outperforms the SJ
method in almost all cases. The evidence suggests that discarding the
normal reference rules,
widely employed by most plug-in rules, can significantly improve the
performance of the plug-in
methods.

The multi-modal test cases 12 through 16 in Table \ref{ISJ table} and
Figure \ref{fig:SJ}
demonstrate that the new bandwidth selection procedure passes the
\textit{bi-modality test} \cite{devroye97}, which
consists of testing the performance of a bandwidth selection procedure
using a bimodal
target density, with the two modes at some distance from each other. It has
been demonstrated in \cite{devroye97} that, by separating the modes of
the target
density enough, existing plug-in selection procedures can be made to
perform arbitrarily poorly due to the adverse effects of the normal
reference rules. The proposed plug-in method in Algorithm \ref{my
plug-in} performs much better than existing plug-in rules, because it
uses the theoretical ideas
developed in \cite{sheather91},
except for the detrimental normal reference rules. A Matlab
implementation of Algorithm \ref{my plug-in} is freely
available from \cite{Matlab}, and includes other examples of improved
performance.

\begin{figure}

\includegraphics{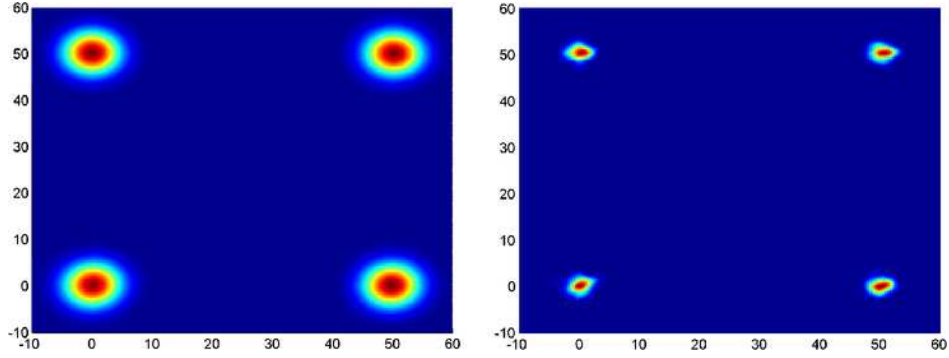}

\caption{Right panel: plug-in rule with normal reference
rule; left panel: the Improved Sheather--Jones method; the normal
reference rule causes significant over-smoothing.}
\label{fig:2dimSJ_oversmooth}
\end{figure}

Algorithm \ref{my plug-in} can be extended to bandwidth selection in
higher dimensions. For completeness we describe the two-dimensional
version of the algorithm in Appendix \ref{appendix:E}. The advantages
of discarding the normal reference rules
persist in the two-dimensional case. In other words, the good
performance of the proposed method in two dimensions
is similar to that observed in the univariate case. For example,
Figure \ref{fig:2dimSJ_oversmooth} shows the superior
%
%
performance of the ISJ method compared to a plug-in
approach using the normal reference rule \cite{wand94,WandJones}, and with
kernels assumed to have a diagonal covariance matrix
with a single smoothing parameter: $\Sigma=t I$. We estimate
the bivariate density,
$\frac{1}{4}\sum_{k=1}^4\mathsf{N}(\bolds\mu_k,I)$, from a
sample of size $N=400$, where
\[
\bolds\mu_1=(0,0),\qquad \bolds\mu_2=(0,50),\qquad
\bolds\mu_3=(50,0),\qquad
\bolds\mu_4=(50,50).
\]
Note that using a plug-in rule with a
normal reference rule causes significant over-smoothing.
The integrated squared error for the ISJ method is 10
times smaller than the corresponding error for the
plug-in rule that uses a normal reference rule \cite{wand94,WandJones}.

\subsection{Bandwidth selection for the diffusion estimator}\label{sec52}

We now discuss the bandwidth choice for the diffusion estimator
(\ref{diff est}). In the following argument we assume that $f$ is as
many times continuously differentiable as needed.
Computation of $t^*$ in (\ref{diffusion_t*}) requires an estimate of
$\Vert Lf\Vert^2$ and $\mathbb E_f[\sigma^{-1}(X)]$. We estimate
$\mathbb
E_f[\sigma
^{-1}(X)]$ via the unbiased estimator
$ \frac{1}{N}\sum_{i=1}^N \sigma^{-1}(X_i)$. The identity
$\Vert Lf\Vert^2=\mathbb E_fL^*Lf(X)$ suggests two possible plug-in
estimators. The
first one is
%
%
\begin{eqnarray}
\widehat{\mathbb E_fL^*Lf}:\!&=&\frac{1}{N}\sum_{j=1}^NL^*Lg(x;t_2)
\Big|_{x=X_j}\nonumber\\[-8pt]\\[-8pt]
&=&\frac{1}{N^2}\sum_{i=1}^N\sum_{j=1}^NL^*L\kappa(x,X_i;t_2)
\Big|_{x=X_j},\nonumber
\end{eqnarray}
where $g(x;t_2)$ is the diffusion estimator (\ref{diff est}) evaluated
at $t_2$, and $\mathscr{X}\equiv\mathbb{R}$. The second estimator is
%
%
\begin{eqnarray}
\label{second stage}
\widehat{\Vert Lf\Vert}{}^2:\!&=&\Vert Lg(\cdot;t_2)\Vert^2\nonumber\\
&=& \biggl\Vert\frac
{\partial
g}{\partial t}(\cdot;t_2)\biggr\Vert^2\nonumber\\[-8pt]\\[-8pt]
&=&\frac{1}{N^2}\sum_{i=1}^N\sum_{j=1}^N\int_\mathbb{R} \,\frac
{\partial\kappa}{\partial t}(x,X_i;t_2)\frac{\partial\kappa
}{\partial t}(x,X_j;t_2) \,dx\nonumber\\
&=&\frac{1}{N^2}\sum_{i=1}^N\sum_{j=1}^NL^*L\kappa(x,X_i;2 t_2)
\Big|_{x=X_j},\nonumber
\end{eqnarray}
where the last line is a simplification that follows from the
Chapman--Kolmogorov equation (\ref{chapman}).
The optimal $t^*_2$ is derived in the same way that
$_*t_2$ is derived for the Gaussian kernel density
estimator. That is, $t^*_2$ is such that both
estimators $\widehat{\mathbb E_fL^*Lf}$ and $\widehat{\Vert Lf\Vert}{}^2$
have the
same asymptotic mean square error.
This leads to the following proposition.
\begin{proposition}
The estimators
$\widehat{\mathbb E_fL^*Lf}$ and $\widehat{\Vert Lf\Vert}{}^2$
have the same asymptotic mean square error when
%
%
\begin{equation}
\label{diffusion second stage}
t^*_2= \biggl(\frac{8+\sqrt{2}}{24} \frac{-3\sqrt{2} \mathbb E_f[\sigma
^{-1}(X)]}{8\sqrt{\pi} N \mathbb E_f[L^*L^2f(X)]} \biggr)^{2/7}.
\end{equation}
\end{proposition}
\begin{pf}
Although the relevant calculations are lengthier, the arguments here
are exactly the same as the ones used in Proposition \ref{prop1}. In particular,
we have the same assumptions on $t$ about its dependence on $N$. For
simplicity of
notation, the operators $L^*$ and $L$ are here assumed
to apply to the first
argument of the kernel~$\kappa$:
\begin{eqnarray*}
&&\mathbb E_f [\widehat{\mathbb E_fL^*Lf} ]\\
&&\qquad=
\mathbb E_f\frac{1}{N^2}\sum_{i=1}^N\sum_{j=1}^NL^*L\kappa(x,X_i;t)
\Big|_{x=X_j}\\
&&\qquad=\frac{1}{N}\int f(x)L^*L\kappa(x,X_i;t) \Big|_{X_i=x} \,dx\\
&&\qquad\quad{} +\frac
{N-1}{N}\iint f(y)f(x)L^*L\kappa(x,y;t)\,dy\,dx\\
&&\qquad=\frac{3\sqrt{2} \mathbb E_f[\sigma^{-1}(X)]}{8\sqrt{\pi} t^{5/2}
N}+O(N^{-1}t^{-3/2})\\
&&\qquad\quad{}+\iint f(y)f(x)L^*L\kappa(x,y;t)\,dy\,dx+O(N^{-1})\\
&&\qquad=\frac{3\sqrt{2} \mathbb E_f[\sigma^{-1}(X)]}{8\sqrt{\pi} t^{5/2}N
}+\int f(y) \int L^*Lf(x)\kappa(x,y;t) \,dx \,dy\\
&&\qquad\quad{} + O\bigl(N^{-1}(1+t^{-3/2})\bigr)\\
&&\qquad=\frac{3\sqrt{2} \mathbb E_f[\sigma^{-1}(X)]}{8\sqrt{\pi} t^{5/2}N
}+\Vert Lf\Vert^2+t\int f(y) L^*L^2f(y) \,dy\\
&&\qquad\quad{}+ O\bigl(N^{-1}(1+t^{-3/2})+t^2\bigr),
\end{eqnarray*}
where we have used a consequence of Lemma \ref{lem:kernel},
\[
\int f(x) L^*L\kappa(x,X_i;t) \Big|_{X_i=x} \,dx \sim\frac{3\sqrt{2}
\mathbb E
_f[\sigma^{-1}(X)]}{8\sqrt{\pi} t^{5/2} },\qquad
t\downarrow0,
\]
and a consequence of the detailed balance equation
(\ref{detailed balance}),
\begin{eqnarray*}
\int L^*Lf(x)\kappa(x,y;t) \,dx&=&\int\frac{p(x)L^*Lf(x)}{p(y)}\kappa
(y,x;t) \,dx\\
&=&L^*Lf(y)+t L^*L^*Lf(y)+O(t^2).
\end{eqnarray*}
Therefore, the squared bias has asymptotic behavior ($N\rightarrow
\infty$)
\[
(\mathbb E_f [\widehat{\mathbb E_fL^*Lf} ]-\Vert Lf\Vert^2 )^2\sim\biggl(\frac
{3\sqrt{2}
\mathbb E_f[\sigma^{-1}(X)]}{8\sqrt{\pi} t^{5/2}N }+t\int f(y) L^*L^2f(y)
\,dy \biggr)^2.
\]
Since\vspace*{1pt} estimator $\widehat{\Vert Lf\Vert}{}^2$ equals
$\widehat{\mathbb E_fL^*Lf}$ when the latter is evaluated at $2 t_2$, the
asymptotic squared bias
of $\widehat{\Vert Lf\Vert}{}^2$ follows immediately, and we simply
repeat the
arguments in the proof of Proposition \ref{prop1} to obtain the desired $t^*_2$.
\end{pf}



Note that $t^*_2$ has the same rate of convergence to $0$
as ${_*}t_2$ in (\ref{_*t_2}). In fact, since
the Gaussian
kernel density estimator is a special case of the diffusion
estimator (\ref{diff est}) when $p(x)=a(x)=1$, the plug-in estimator
(\ref{second stage}) for the
estimation of $\Vert Lf\Vert^2$ reduces to the plug-in
estimator for the estimation of
$\frac{1}{4}\Vert f''\Vert^2$. In addition, when $p(x)=a(x)=1$, the $t_2^*$
in (\ref{diffusion second stage}) and ${_*}t_2$ in (\ref{_*t_2}) are
identical.
We thus
suggest the following bandwidth selection and estimation procedure for
the diffusion estimator~(\ref{diff est}).
\begin{algorithm}
\label{diffusion procedure}
%
\begin{enumerate}
\item Given the data $X_1,\ldots,X_N$, run Algorithm \ref{my plug-in}
to obtain
the Gaussian kernel density estimator (\ref{KDE}) evaluated at $_*\hat
t$ and the optimal bandwidth $\sqrt{{_*}\hat t_2}$ for the estimation
of $\Vert f''\Vert^2$.
This is the pilot estimation step.
\item Let $p(x)$ be the Gaussian kernel density estimator from step 1, and
let $a(x)=p^\alpha(x)$ for some $\alpha\in[0,1]$.
\item Estimate $\Vert Lf\Vert^2$ via the plug-in estimator (\ref{second
stage}) using $\hat t_2^*={_*}\hat t_2$, where ${_*}\hat t_2$ is
computed in step 1.
\item Substitute the estimate of $\Vert Lf\Vert^2$ into
(\ref{diffusion_t*}) to obtain an estimate for $t^*$.
\item Deliver the diffusion estimator (\ref{diff est})
evaluated at $\hat t^*$ as the final density estimate.
\end{enumerate}
\end{algorithm}

The bandwidth selection rule that we use for the diffusion estimator in
Algorithm \ref{diffusion procedure}
is a single stage direct plug-in
bandwidth selector, where the bandwidth $t^*_2$ for the estimation of
the functional
$\Vert Lf\Vert^2$ is approximated by $_*\hat t_2$ (which is computed in
Algorithm \ref{my plug-in}), instead of being derived from a normal
reference rule.
In the next section, we illustrate
the performance of Algorithm \ref{diffusion procedure} using some well-known
test cases for density estimation.
%
%
\begin{remark}[(Random variable generation)]
For applications of kernel density estimation, such as the smoothed
bootstrap, efficient random variable generation from the diffusion estimator
(\ref{diff est}) is
accomplished via the Euler method as applied to the stochastic
differential equation
(\ref{Ito}) (see \cite{Kloeden}).
\begin{algorithm}
\begin{enumerate}
\item Subdivide the interval $[0,\hat t^*]$ into $n$ equal intervals of
length $\delta t=\hat t^*/n$ for some large $n$.
\item Generate a random integer $I$ from $1$ to $N$ uniformly.
\item For $i=1,\ldots,n$, repeat
\[
Y_i=Y_{i-1}+\mu(Y_{i-1}) \delta t+\sigma(Y_{i-1})\sqrt{\delta t}
Z_i,
\]
where $Z_1,\ldots,Z_n\sim_{\mathrm{i.i.d.}}\mathsf{N}(0,1)$, and $Y_0=X_I$.
\item Output $Y_n$ as a random variable with approximate density (\ref
{diff est}).
\end{enumerate}
\end{algorithm}

Note that since we are only interested in the
approximation of the statistical properties of $Y_n$,
there are no gains to be had from using the more complex
Milstein stochastic integration procedure \cite{Kloeden}.
\end{remark}

\section{Numerical experiments}\label{sec6}

In this section, we provide a simulation study of the diffusion
estimator. In implementing Algorithm
\ref{diffusion procedure}, there are a number of issues to consider.
First, the numerical solution of the PDE (\ref{diffusion}) is a straightforward
application of either finite difference or spectral methods
\cite{LarssonThomee}. A Matlab implementation using finite differences
and the
stiff ODE solver \texttt{ode15s.m} is available from the first author
upon request. Second, we compute
$\Vert Lg(\cdot;\hat t_2^*)\Vert^2$ in
Algorithm \ref{diffusion procedure} using the approximation
\[
\Vert Lg(\cdot; t)\Vert^2= \biggl\Vert\frac{\partial g}{\partial t}(\cdot;t)
\biggr\Vert^2\approx\Vert g(\cdot;t+\varepsilon)-g(\cdot;t)
\Vert^2/\varepsilon^2,\qquad
\varepsilon\ll1,
\]
where $g(\cdot;t)$ and $g(\cdot;t+\varepsilon)$ are the successive
output of the
numerical integration routine (\texttt{ode15s.m} in our case).
Finally, we selected $\alpha=1$ or $a(x)=p(x)$ in Algorithm \ref
{diffusion procedure} without using any clipping of the pilot estimate.
For a small simulation study with $\alpha=0$,
see \cite{preprint}.

We would like to point out that simulation studies of existing
variable-location scale
estimators \cite{jones94,Samiuddin,TerrelScott} are implemented
assuming that the target p.d.f. $f$ and any functionals of $f$ are known
\textit{exactly} and no pilot estimation step is employed. In addition,
in these simulation studies the bandwidth is chosen so that it is the
global minimizer of the exact MISE. Since in practical applications the
MISE and all functionals of $f$ are not available, but have to be
estimated, we proceed differently in our simulation study.
We compare the estimator of Algorithm \ref{diffusion procedure} with
the Abramson's popular adaptive kernel density estimator
\cite{abramson}. The parameters $_*t$ and
$_*t_2$ of the diffusion estimator are estimated using the new
bandwidth selection procedure in Algorithm \ref{my plug-in}. The
implementation of Abramson's estimator in the
\textsf{Stata} language is given in \cite{stata}.
Briefly, the estimator is given by
\[
\hat f_\mathsf{A}(x)=\frac{1}{N \sqrt{t}\lambda_i}\sum_{i=1}^N
\phi\biggl(\frac{x-X_i}{\sqrt{t}\lambda_i} \biggr),
\]
where $\lambda_i^2=G/\hat f(X_i;t_p)$, $G= (\prod_{i=1}^N \hat
f(X_i;t_p) )^{1/N}$, and the bandwidths $\sqrt{t}$ and $\sqrt{t_p}$
are computed using Least Squares Cross Validation (LSCV)
\cite{Loader}.

Our criterion for the comparison is the numerical approximation to
\[
\operatorname{Ratio}=\frac{\Vert g(\cdot;\hat t^*)-f\Vert^2}{\Vert\hat
f_\mathsf
{A}-f\Vert^2},
\]
that is, the ratio of the integrated squared error of the diffusion
estimator to the integrated squared error of the alternative kernel
density estimator.

%
\begin{table}[b]
\tabcolsep=0pt
\caption{Results over 10 independent simulation experiments. In all
cases the domain was assumed to be $\mathbb{R}$}\label{table2}
\begin{tabular*}{\tablewidth}{@{\extracolsep{\fill}}lcccc@{}}
\hline
\textbf{Case} & \textbf{Target density} $\bolds{f(x)}$ & $\bolds{N}$ & \textbf{Ratio I} & \textbf{Ratio II} \\
\hline
1&$\frac{1}{2}\mathsf{N} (0, (\frac{1}{10} )^2 )+\frac{1}{2}\mathsf
{N} (5,1 )$ &
$10^3$ & 0.9 & 0.82\\
& & $10^5$ & 0.23 & 0.48\\
[4pt]
2&$\frac{1}{2}\mathsf{N}(0,1)+\sum_{k=0}^4\frac{1}{10}\mathsf{N}
(\frac{k}{2}-1, (\frac{1}{10} )^{2} )$ & $10^3$
& 0.65 & 0.99\\
& & $3\times10^5$ & 0.11 & 0.51\\
[4pt]
3&$\sum_{k=0}^7\frac{1}{8}\mathsf{N} (3 ( (\frac{2}{3} )^k-1 ),
(\frac{2}{3} )^{2k} )$ &
$10^3$ & 1.05 & 0.75\\
& & $10^5$ & 0.15 & 0.45\\
[4pt]
4&$\frac{49}{100}\mathsf{N} (-1, (\frac{2}{3} )^2 )+\frac
{49}{100}\mathsf{N} (1, (\frac{2}{3} )^2 )+\frac{1}{350}\sum
_{k=0}^6\mathsf{N} (\frac{k-3}{2}, (\frac{1}{100} )^2 )$ &
$10^3$ & 0.94 & 0.63\\
& & $10^5$ & 0.46 & 0.76\\
[4pt]
5&$\frac{2}{7}\sum_{k=0}^2\mathsf{N} (\frac
{12k-15}{7}, (\frac{2}{7} )^2 )+\frac{1}{21}\sum_{k=8}^{10} \mathsf
{N} (\frac{2k}{7}, (\frac{1}{21} )^2 ) $ &
$10^3$ & 0.54 & 2.24\\
& & $10^5$ & 0.12 & 0.84\\
[4pt]
6&$\frac{46}{100}\sum_{k=0}^1 \mathsf{N} (2k-1, (\frac
{2}{3} )^2 )+
\sum_{k=1}^3\frac{1}{300}\mathsf{N} (-\frac{k}{2}, (\frac{1}{100}
)^2 )$ & $10^4$ &
0.83 & 0.93\\[4pt]
& \multicolumn{1}{c}{\hspace*{-105pt}\hspace*{48pt}$\,+\sum_{k=1}^3\frac{7}{300}\mathsf{N} (\frac{k}{2}, (\frac
{7}{100} )^2 )$} & $10^5$ & 0.55 & 0.68\\
[4pt]
7&$ \frac{1}{2}\mathsf{N}(-2,\frac{1}{4})+\frac{1}{2}\mathsf
{N}(2,\frac{1}{4})$ & $10^3$ &
0.51 & 0.51\\
& & $10^5$ & 0.41 & 0.89\\
[4pt]
8&$\frac{3}{4}\mathsf{N}(0,1)+\frac{1}{4}\mathsf{N} (\frac{3}{2},
(\frac{1}{3} )^2 )$ &
$10^3$ & 0.59 & 0.53\\
& & $10^6$ & 0.79 & 1.01\\
[4pt]
9&Log-Normal with $\mu=0$ and $\sigma=1$ &
$10^3$ & 0.17 & 0.85\\
& & $10^5$ & 0.12 & 0.51\\
[4pt]
10&$ \frac{1}{2}\mathsf{N}(0,1)+\sum_{k=-2}^2\frac
{2^{1-k}}{31}\mathsf{N} (k+\frac{1}{2}, (\frac{2^{-k}}{10} )^2 ) $ &
$10^3$ & 0.88 & 0.98\\
& & $10^4$ & 0.30 & 0.85\\
\hline
\end{tabular*}
\end{table}


%
%

Table \ref{table2}, column 4 (ratio I) shows
the average results over 10 independent trials for a number of different
test cases. The second column displays the target density and the
third column shows the sample size used for the experiments. In the
table $\mathsf{N}(\mu,\sigma^2)$, denotes a Gaussian density with
mean $\mu$ and variance $\sigma^2$. Most
test problems are taken from \cite{MarronWand}.
For each test case, we conducted a simulation run with both a
relatively small sample size and a relatively large sample size
wherever possible.
The table shows that, unlike the standard
variable location-scale estimators \cite{jones94,TerrelScott}, the
diffusion estimator does not require any clipping
procedures in order to retain its good performance for large
sample sizes.

Next, we compare the practical performance of the proposed diffusion
estimator with the performance of higher-order kernel estimators. We
consider the sinc kernel estimator defined as
\[
\hat f_\mathsf{sinc}(x)=\frac{1}{N}\sum_{i=1}^N \frac{1}{\sqrt
{t}}K \biggl(\frac{x-X_i}{\sqrt{t}} \biggr),\qquad K(x)=\frac{\sin(x)}{\pi x},
\]
where again $\sqrt{t}$ is selected using LSCV. Table \ref{table2}, column 5
(ratio II) shows that the results are broadly similar and our method is
favored in all cases except test case~5.
Higher-order kernels do not
yield proper density estimators, because
the kernels take on negative values. Thus, an important advantage of
our method
and all second order kernel methods is that they provide nonnegative
density estimators.
As pointed out in \cite{WandJones}, the good asymptotic performance of
higher-order kernels
is not guaranteed to carry over to finite sample sizes in practice. Our
results confirm this observation.

%
\begin{table}
\caption{Practical performance of the boundary bias correction of the
diffusion estimator for the test cases: (1) exponential distribution
with mean equal to unity; (2) test cases 1 through 8, truncated to the
interval $(-\infty,0]$}
\label{table3}
\begin{tabular*}{\tablewidth}{@{\extracolsep{\fill}}lccccccccc@{}}
\hline
\textbf{Test case} & $\bolds{\mathsf{Exp}(1)}$ & \textbf{1}
& \textbf{2} & \textbf{3} & \textbf{4} & \textbf{5} & \textbf{6} & \textbf{7} & \textbf{8} \\
\hline
Ratio & 0.52 & 0.38 &0.74 & 0.25 & 0.70 & 0.38 & 0.74 & 0.56 & 0.46 \\
\hline
\end{tabular*}
\end{table}

In addition, we make a comparison with
the novel polynomial boundary correction method of Hall and Park
\cite{hallpark}. The results are given
in Table \ref{table3}, where we use some of the test cases defined in
Table \ref{ISJ table}, truncated to the interval $(-\infty,0]$. Table~\ref{table3} shows that
for finite sample sizes the practical performance of our approach is
competitive.
We now give the implementation details.
Let $\beta$ be the
point of truncation from above, which is assumed to be known in advance.
Then, the Hall and Park estimator is
%
%
\begin{equation}
\label{HP}\qquad
\hat f_\alpha(x;t)=\frac{1}{N \int_{-\infty}^\beta\phi((x-y)/h)\,d y}
\sum_{i=1}^N\phi\biggl(\frac{x-X_i+\alpha(x)}{\sqrt{t}} \biggr),\qquad x\leq\beta,
\end{equation}
where $ \alpha(x) =t \frac{\hat f'_0(x)}{\hat f_0(x)}\rho(\frac
{x-a}{h} );$ $\hat f_0(x)$
is equivalent to $\hat f_\alpha(x)$ when $\alpha(x)\equiv0$, and
$\hat f'_0(x)$ is an estimator of $f'(x)$; $\rho(u)=\frac{1}{\phi
(u)}\int_{-\infty}^u v \phi(v)\,dv$.
We use LSCV
to select a suitable bandwidth $\sqrt{t}$.
The denominator in (\ref{HP})
adjusts for the deficit of probability mass in the neighborhood
of the end-point, but note that theoretically
(\ref{HP}) does not integrate to unity and therefore random variable
generation from (\ref{HP}) is not straightforward. In addition,
our estimator more easily handles the case with two end-points. On the
positive side, Hall and Park
\cite{hallpark} note that their estimator preserves
positivity and has excellent asymptotic properties, which
is an advantage over many other boundary kernels.

Finally, we give a two-dimensional density estimation example, which to
the best of our knowledge cannot be handled satisfactorily by existing
methods \cite{HallPark02,Alberts} due to the boundary bias effects.
%
%
\begin{figure}[b]

\includegraphics{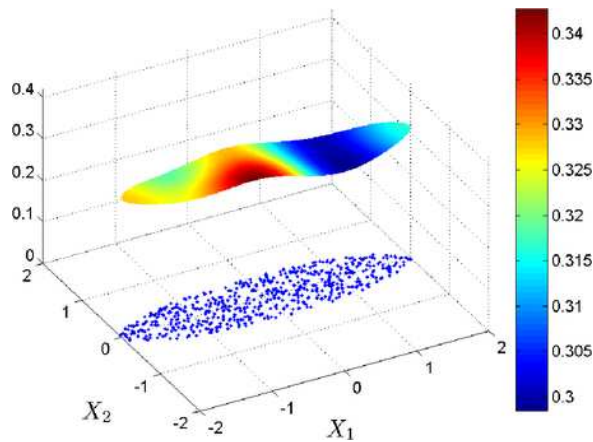}

\caption{A two-dimensional example with 600 points generated uniformly
within an ellipse.}
\label{fig:2d}
\end{figure}
The two-dimensional version
of equation (\ref{CH}) is
\begin{eqnarray*}
\frac{\partial\hat f}{\partial t}(\mathbf{x};t)&=&\frac{1}{2} \biggl(\frac
{\partial^2\hat f}{\partial x_1^2}(\mathbf{x};t)+\frac{\partial
^2\hat f}{\partial x_2^2}(\mathbf{x};t) \biggr)\qquad \forall t>0, \mathbf{x}\in
\mathscr{X}, \\
\hat f(\mathbf{x};0)&=&\Delta(\mathbf{x}),\\
\mathbf{n}\cdot\nabla\hat f(\mathbf{x};t) &=&0\qquad \forall t>0,
\end{eqnarray*}
where $\mathbf{x}=(x_1,x_2)$ belongs to the set
$\mathscr{X}\subseteq\mathbb{R}^2$, the initial condition $\Delta
(\mathbf{x})$ is the
empirical density of the data, and in the Neumann boundary condition
$\mathbf{n}$ denotes the unit outward normal to the boundary $\partial
\mathscr{X}$ at $\mathbf{x}$. The particular example which we
consider is the density estimation
of $600$ uniformly distributed points on the domain $\mathscr{X}=\{
\mathbf{x}\dvtx
x_1^2+(4x_2)^2\leq4\}$. We assume that the domain of the data
$\mathscr{X}$ is known prior to the estimation. Figure
\ref{fig:2d} shows $\hat f(\mathbf{x};\hat t^*)$ on
$\mathscr{X}=\{\mathbf{x}\dvtx x_1^2+(4x_2)^2\leq4\}$, that is, it
shows the numerical solution of the two-dimensional PDE
at time $\hat t^*=0.13$ on the set $\mathscr{X}$. The
bandwidth was determined using the bandwidth selection
procedure described in Appendix \ref{appendix:E}. We emphasize the
satisfactory way in which the p.d.f. $\hat f(\mathbf{x};\hat t^*)$
handles any boundary bias problems.
It appears that currently existing methods
\cite{HallPark02,foster,Alberts,Karunamuni} cannot handle such
two-dimensional (boundary) density estimation problems either because
the geometry of the set $\mathscr{X}$ is too complex, or because the
resulting estimator is not a bona-fide p.d.f.

\section{Conclusions and future research}\label{sec7}

We have presented a new kernel density estimator based
on a linear diffusion process. The key idea is to
construct an adaptive kernel by considering the most
general linear diffusion with its
stationary density equal to a pilot density estimate.
The resulting diffusion estimator unifies many of
the existing ideas about adaptive smoothing. In
addition, the estimator is consistent at boundaries.
Numerical experiments suggest
good practical performance. As future research,
the proposed estimator can
be extended in a number of ways.
First, we can construct kernel density estimators
based
on L\'evy processes, which will have the diffusion
estimator as a special case. The kernels
constructed via
a L\'evy process could be tailored for data for
which
smoothing with the Gaussian kernel density estimator or diffusion
estimator is
not
optimal. Such cases arise when the data is a sample from a heavy-tailed
distribution. Second, more subtle and interesting smoothing models can
be constructed by considering nonlinear parabolic
PDEs.
One such candidate is the quasilinear parabolic PDE
with diffusivity that depends on the density
exponentially:
\[
\frac{\partial}{\partial t}g(x;t)=\frac{\partial}{\partial x} \biggl(\e
^{-\alpha g(x;t)}\,\frac{\partial}{\partial x}g(x;t) \biggr),\qquad \alpha>0.
\]
Another viable model is the semilinear parabolic PDE
\[
\frac{\partial}{\partial t} \bigl(\e^{u(x;t)} \bigr)=\frac{1}{2}\,\frac
{\partial^2}{\partial x^2}u(x;t),
\]
where $u(x;t)=\log(g(x;t))$ is the logarithm of the
density estimator. The Cauchy density
$\frac{t}{\pi(x^2+t^2)}$ is a particular solution and thus
the model could be useful for smoothing heavy-tailed data.
All such nonlinear models will provide adaptive smoothing without
the need for a pilot run, but at the cost of increased
model complexity.

\begin{appendix}
\section{Gaussian kernel density estimator properties}
\label{append:A}

In this appendix, we present the technical
details for the proofs of the properties of the diffusion estimator. In
addition, we include a description of our plug-in rule in two dimensions.

We use \mbox{$\Vert\cdot\Vert$} to denote the Euclidean norm on $\mathbb{R}$.
\begin{theorem} Let $t=t_N$ be such that $\lim_{N\rightarrow\infty
}t_N=0$ and $\lim_{N\rightarrow\infty} N\sqrt{t_N}=\infty$. Assume
that $f''$ is a continuous square-integrable function. The integrated
squared bias and integrated variance of the Gaussian kernel density
estimator (\ref{KDE}) have asymptotic behavior
%
%
\begin{equation}
\label{kde_bias}\Vert\mathbb E_f[\hat f( \cdot;t)]-f\Vert^2=\tfrac
{1}{4}t^2\Vert f''\Vert^2+o(t^2),\qquad N\rightarrow\infty,
\end{equation}
and
%
%
\begin{equation}
\label{kde_var}
\int\operatorname{Var}_f[\hat f(x;t)] \,dx=\frac{1}{2N\sqrt{\pi
t}}+o\bigl(\bigl(N\sqrt{t}\bigr)^{-1}\bigr),\qquad N\rightarrow\infty,
\end{equation}
respectively. The first-order asymptotic approximation of MISE, denoted
AMISE, is thus given by
%
%
\begin{equation}\label{AMISE}
\operatorname{AMISE}\{\hat f\}(t)=\frac
{1}{4}t^2\Vert f''\Vert^2+\frac{1}{2N\sqrt{\pi t}}.
\end{equation}
The asymptotically optimal value of $t$ is the minimizer of the AMISE
%
%
\begin{equation}
\label{Gaussian kde bandwidth} _*t= \biggl(\frac{1}{2N\sqrt{\pi}
\Vert f''\Vert^2} \biggr)^{2/5},
\end{equation}
giving the minimum value
%
%
\begin{equation}\label{rate}
\operatorname{AMISE}\{\hat f\}(_*t)=N^{-4/5}
\frac{5\Vert f''\Vert^{2/5}}{4^{7/5}\pi^{2/5}} .
\end{equation}
\end{theorem}

For a simple proof, see \cite{WandJones}.

\section{\texorpdfstring{Proof of Lemma
\protect\lowercase{\ref{lem:kernel}}}{Proof of Lemma 1}}
\label{appendix:B}

We seek to establish the behavior of the solution of (\ref{BE}) and
(\ref{FE}) as $t\downarrow0$.
We use the Wentzel--Kramers--Brillouin--Jeffreys (WKBJ) method described
in \cite{azencott,Cohen,kannai,molchanov}. In the WKBJ method, we look
for an asymptotic expansion of the form
%
%
\begin{equation}
\label{asymptotic expansion}
\kappa(x,y;t)\sim\e^{-{1}/({2t})s^2(x,y)}\sum_{m=0}^\infty
t^{m-1/2}C_m(x,y),\qquad t\downarrow0,
\end{equation}
where $\{C_m(x,y)\}$ and $s(x,y)$ are unknown functions.
To determine $s(x,y)$ and $\{C_m(x,y)\}$, we
substitute the expansion into (\ref{FE}) and, after
canceling the exponential term, equate coefficients of
like powers of $t$. This matching of the powers of $t$ leads to
solvable ODEs, which determine the unknown functions.
Eliminating the leading order $O(t^{-5/2})$ term gives the ODE for $s$
%
%
\begin{equation}
\label{ode1}
a(x) \biggl[\frac{\partial}{\partial x}s(x,y) \biggr]^2-p(x)=0.
\end{equation}
Setting the next highest order $O(t^{-3/2})$ term in the expansion
to zero gives the ODE
%
%
\begin{eqnarray}
\label{ode2}\qquad
0 &=& 2 a ( x ) s (
x,y ) \,{\frac{\partial s }{\partial x}} \,{\frac{d p}{dx}} p
( x ) C_0 ( x,y )-2 a ( x ) s
( x,y ) \,{\frac{\partial s }{\partial x}} p ^2 ( x ) \,{\frac{
\partial C_0}{\partial x}}\nonumber\\
&&{} + p^3 ( x ) C_0 ( x,y ) + s ^2 ( x,y ) p ^3 ( x ) C_1
( x,y )\nonumber\\
&&{}  - {\frac{d a}{dx}}
p^2 ( x ) s ( x,y )\,
{\frac{\partial s}{\partial x}} C_0
( x,y ) \\
&&{}+a ( x ) s^2 ( x,y )
\biggl( {\frac{\partial s }{\partial x}} \biggr) ^{2} p^2 ( x ) C_1
( x,y ) -a ( x ) \biggl( {\frac{\partial s}{
\partial x}} \biggr) ^{2} p^2 ( x
) C_0 ( x,y )\nonumber\\
&&{} -a ( x ) s
( x,y )\, {\frac{\partial^{2} s}{\partial{x}^{2}}} p^2 ( x ) C_0
( x,y ).\nonumber
\end{eqnarray}
To determine a unique solution to (\ref{ode1}), we impose the
condition $s(x,x)=0$, which is necessary, but not sufficient, to ensure
that $\lim_{t\downarrow0}\kappa(x,y;t)=\delta(x-y)$. This gives the solution
\[
s(x,y)=\int_y^x\sqrt{\frac{p(s)}{a(s)}}\,ds.
\]
Substituting this solution into (\ref{ode2})
and simplifying gives an equation without $C_1(x,y)$,
%
%
\begin{equation}
C_0(x,y) p(x) \,\frac{da}{dx}+4 a(x) p(x) \,\frac{\partial C_0}{\partial
x}-3 C_0(x,y) \,\frac{dp}{dx} a(x)=0,
\end{equation}
whence we have the general solution
$
C_0(x,y)=h(y) p^{3/4}(x) a^{-1/4}(x)
$
for some as yet unknown function of $y$, $h(y)$. To determine $h(y)$,
we require that
the kernel $\widetilde\kappa(x,y;t)$ satisfies the detailed balance equation
(\ref{detailed balance}). This ensures that $\widetilde\kappa
(x,y;t)$ also satisfies (\ref{BE}). It follows that
$C_0(x,y)$ has to satisfy $p(y)C_0(x,y)=p(x)C_0(y,x)$, which after
rearranging gives
\[
h(x)(a(x)p(x))^{1/4}=h(y)(a(y)p(y))^{1/4}.
\]
A separation of variables argument now gives
$h(y)(a(y)p(y))^{1/4}=\mathrm{const}.$, and
hence
\[
C_0(x,y)=\mathrm{const.} (a(y)p(y))^{-1/4} p^{3/4}(x) a^{-1/4}(x).
\]
We still need to determine the arbitrary constant. The constant is
chosen so that
\[
\lim_{t\downarrow0} \int_{-\infty}^\infty\widetilde\kappa(x,y;t)
\,dx=1,
\]
which ensures that $\lim_{t\downarrow0}\widetilde\kappa
(x,y;t)=\delta(x-y)$.
This final condition yields
\[
C_0(x,y)=\frac{p(x)}{\sqrt{2\pi}(a(y)p(y)a(x)p(x))^{1/4}},
\]
and hence
\[
\widetilde\kappa(x,y;t)=\frac{p(x)}{\sqrt{2\pi t}
[p(x)a(x)a(y)p(y)]^{1/4}}\exp\Biggl\{-\frac{1}{2t} \Biggl[\int_y^x\sqrt{\frac
{p(s)}{a(s)}} \,ds \Biggr]^2 \Biggr\}.
\]
%
%
%
\begin{remark}
Matching higher powers of $t$ gives first order linear ODEs for the
rest of the unknown functions $\{C_m(x,y),m\geq1\}$. The ODE for each
$C_m(x,y), m=1,2,3,\ldots$ is
\[
as'(C_m/p)'+ \biggl(\frac{(as')'}{2p}+(m-1/2)
\biggr)C_m=\bigl(a(C_{m-1}/p)'\bigr)',\qquad
C_m(y,y)=0,
\]
where all derivatives apply to the variable $x$ and $y$ is treated as a
constant. Thus, in principle, all functions $\{C_m(x,y)\}$ can be
uniquely determined.
\end{remark}

It can be shown (see \cite{Cohen}) that the expansion
(\ref{asymptotic expansion}) is valid under the conditions that $a,p$
and all their derivatives are bounded from above, and $p(x)\geq p_0>0$,
\mbox{$a(x)\geq a_0>0$}. Here, we only establish the validity of the leading
order approximation
$\widetilde\kappa$ under the milder conditions (\ref{big
assumption}). We do not attempt to prove the validity of the higher
order terms in
(\ref{asymptotic expansion}) under the weaker conditions. The proof of
the following lemma uses arguments similar to the ones given in \cite{Cohen}.
\begin{lemma}
\label{Joe's lemma}
Let $a(x)$ and $p(x)$ satisfy conditions (\ref{big assumption}).
Then, for all $t\in(0,t_0]$, where
$t_0>0$ is some constant independent of $x$ and $y$, there holds
\[
|\kappa(x,y;t)-\widetilde\kappa(x,y;t)|\leq \operatorname{const.} C_0(x,y)
t^{1/4}\e^{-{s^2(x,y)}/({2t})}\qquad \forall x,y.
\]
\end{lemma}

To prove the lemma, we first begin by proving the following auxiliary results.
\begin{proposition}
\label{prop:bound}
Define
\[
\ell(z)=\ell(z;x,y,t,\tau)=\frac{s^2(x,z)}{2(t-\tau)}+\frac
{s^2(z,y)}{2\tau}.
\]
Then for $\tau\in(0,t)$, we have
\[
\ell(z)\geq\frac{s^2(x,y)}{2t}.
\]
Moreover, there exists
a unique $z_0=z_0(x,y,t,\tau)$ for which
$\ell(z_0)=\frac{s^2(x,y)}{2t}$, and
$\ell(z)$ is increasing for $z>z_0$
and decreasing for $z<z_0$.
\end{proposition}
\begin{pf}
We have
\[
\ell(z)=\frac{1}{2(t-\tau)} \biggl(\int_z^x \sigma^{-1}(s) \,ds \biggr)^2+\frac
{1}{2\tau} \biggl( \int_y^z \sigma^{-1}(s) \,ds \biggr)^2,
\]
and hence
%
%
\begin{equation}
\label{deriv}
\ell'(z)=\frac{-\sigma^{-1}(z)}{t-\tau}\int_z^x\sigma^{-1}(s)
\,ds+\frac{\sigma^{-1}(z)}{\tau}\int_y^z\sigma^{-1}(s) \,ds.
\end{equation}
For $x\not=y$, $\ell'(y)>0, \ell'(x)<0$, and therefore
by the continuity of $\ell'$,
there exists $z_0\in(x,y)\dvtx\ell'(z_0)=0$.
For $x=y$, set $z_0=x$. Setting $z=z_0$ in
(\ref{deriv}),
%
%
\begin{equation}
\label{star}
\frac{1}{t-\tau}\int_{z_0}^x \sigma^{-1}(s)\, ds=\frac{1}{\tau}\int
_y^{z_0}\sigma^{-1}(s)\, ds.
\end{equation}
Therefore,
${\int_{z_0}^x}\sigma^{-1}(s) \,ds=\frac{t-\tau}{\tau}{\int
_y^{z_0}}\sigma^{-1}(s) \,ds$ and adding
${\int_y^{z_0}}\sigma^{-1}(s) \,ds$ to both sides we obtain
\[
\int_y^{x}\sigma^{-1}(s) \,ds= \frac{t}{\tau}\int_y^{z_0}\sigma
^{-1}(s) \,ds,
\]
from which we see that (\ref{star}) is also equal to
$\frac{1}{t}{\int_y^{x}}\sigma^{-1}(s) \,ds$.
Hence, by substitution $\ell(z_0)= \frac{1}{2t} (\int_{y}^{x}\sigma
^{-1}(s) \,ds )^2$,
as required. Finally, note that if
$F(z)=\ell(z)-\frac{t}{2\tau(t-\tau)} (\int_{z_0}^z\sigma^{-1}(s)
\,ds )^2$, then\vspace*{1pt}
$F'(z)=0$ for all $z$. Hence, $F(z)=F(z_0)=\ell(z_0)$ and
%
%
\begin{equation}
\label{fortytwo}
\ell(z)=\ell(z_0)+\frac{t}{2\tau(t-\tau)} \biggl(\int_{z_0}^z\sigma
^{-1}(s) \,ds \biggr)^2.
\end{equation}
As a consequence of Proposition \ref{prop:bound}, we have the
following result.
\end{pf}
\begin{proposition} Assuming $\lim_{z\rightarrow\pm\infty}\int
_{z_0}^z\sigma^{-1}(s)\,ds=\pm\infty$, we have the following equality:
\label{prop:bound2}
\begin{eqnarray*}
&&\int_0^t\sqrt{\int_{-\infty}^\infty\biggl( \frac{\e^{-
{s^2(x,z)}/({2(t-\tau)})}}{\sqrt{t-\tau}} \frac{\e^{-
{s^2(z,y)}/({2\tau})}}{\sqrt{\sigma(z)}\sqrt{\tau}} \biggr)^2 \,dz} \,d\tau
\\
&&\qquad=2\pi^{-1/4}t^{1/4} \Gamma^2(3/4) \e^{-
{s^2(x,y)}/({2t})}\\
&&\qquad= c_2 t^{1/4}\e^{-{s^2(x,y)}/({2t})},
\end{eqnarray*}
where $c_2$ is a constant [indeed $c_2=2\pi^{-1/4}\Gamma^2(3/4)$].
\end{proposition}
\begin{pf}
We have
\begin{eqnarray*}
&&
\int_{-\infty}^\infty\frac{\e^{-{s^2(x,z)}/{(t-\tau
)}}}{t-\tau} \frac{\e^{-{s^2(z,y)}/{\tau}}}{\sigma(z)\tau}
\,dz\\
&&\qquad= \frac{1}{(t-\tau)\tau} \e^{-2\ell(z_0)}\int_{-\infty
}^\infty\sigma^{-1}(z)\e^{-{ (\int_{z_0}^z\sigma^{-1}(s)\,ds
)^2}/({\tau(t-\tau)/t})}\,dz\\
&&\qquad=\frac{1}{\sqrt{t(t-\tau)\tau}} \e^{-2\ell(z_0)}\int_{-\infty
}^\infty\e^{-v^2}\,dv,
\end{eqnarray*}
with the change of variable $ v(z)=\frac{1}{\sqrt{\tau(t-\tau
)/t}}\int_{z_0}^z\sigma^{-1}(s) \,ds$. Then the result follows
from the fact that $\int_0^t(\tau(t-\tau))^{-1/4}\,d\tau=2\pi
^{-1/2}t^{1/2}\Gamma^2(3/4)$.

Given these two auxiliary results, we proceed with the proof of Lemma
\ref{Joe's lemma}. Writing
\[
\kappa^*(x,y;t)=
\frac{\partial}{\partial t}\widetilde\kappa(x,y;t)-L\widetilde
\kappa(x,y;t)=-\frac{\e^{-{s^2(x,y)}/({2t})}}{\sqrt{t}}LC_0(x,y),
\]
we define inductively the following sequence of function $\{\rho_j\}$,
starting with \mbox{$\rho_0=0$}:
\begin{eqnarray}
\rho_{j+1}(x,y;t)&=&-\kappa^*(x,y;t)-\int_0^t\int_{-\infty}^\infty
\kappa^*(x,z;t-\tau)\rho_j(z,y;\tau)\,dz\,d\tau,\nonumber\\
\eqntext{j=1,2,\ldots.}
\end{eqnarray}
Note in particular that $\rho_1=-\kappa^*$.
We will show that there exists a limit of $\{\rho_j\}$.
We begin by proving via induction that for $j\geq1$, $x,y\in\mathbb
{R}$, $t\in(0,t_0]$,
where
\[
t_0=\min\biggl\{ \biggl(\frac{\sqrt{2\pi}}{2c_1c_2} \biggr)^{4/3},1 \biggr\},
\]
there holds
%
%
\begin{equation}
\label{indo}
|\rho_{j+1}(x,y,t)-\rho_j(x,y,t)|\leq\frac
{c_3}{2^j}|LC_0(x,y)|t^{1/4}\e^{-{s^2(x,y)}/({2t})},
\end{equation}
where $c_3=2c_1c_2/\sqrt{2\pi}$.
First, we calculate
for $j=1$
\[
\rho_2(x,y,t)=-\kappa^*(x,y,t)+\int_0^t\int_{-\infty}^\infty
\kappa^*(x,z,t-\tau)\kappa^*(z,y,\tau)\,dz\,d\tau.
\]
Therefore, we have the following bound:
\begin{eqnarray*}
&&
|\rho_2(x,y,t)-\rho_1(x,y,t)|\\
&&\qquad\leq
\int_0^t\int_{-\infty}^\infty|\kappa^*(x,z,t-\tau)\kappa
^*(z,y,\tau)|\,dz\,d\tau\\
&&\qquad= \int_0^t\int_{-\infty}^\infty\frac{\e^{-
{s^2(x,z)}/({2(t-\tau)})}}{\sqrt{t-\tau}}
\frac{\e^{-{s^2(z,y)}/({2\tau})}}{\sqrt{\tau
}}|LC_0(x,z)LC_0(z,y)|\,dz\,d\tau\\
&&\qquad= \int_0^t\int_{-\infty}^\infty\frac{\e^{-
{s^2(x,z)}/({2(t-\tau)})}}{\sqrt{t-\tau}}
\frac{\e^{-{s^2(z,y)}/({2\tau})}}{\sqrt{\sigma(z)\tau}}\\
&&\qquad\quad\hspace*{34.8pt}{}\times\sqrt
{\sigma(z)}|LC_0(x,y)|\frac{|Lq(z)|}{\sqrt{2\pi
}(a(z)p(z))^{1/4}}\,dz\,d\tau\\
&&\qquad = \frac{1}{\sqrt{2\pi}} |LC_0(x,y)|\int_0^t\int_{-\infty
}^\infty\frac{\e^{-{s^2(x,z)}/({2(t-\tau)})}}{\sqrt{t-\tau}}\\
&&\qquad\quad\hspace*{113.6pt}{}\times
\frac{\e^{-{s^2(z,y)}/({2\tau})}}{\sqrt{\sigma(z)\tau}} \frac
{|Lq(z)|}{ q(z) }\,dz\,d\tau\\
&&\qquad\leq\frac{1}{\sqrt{2\pi}} |LC_0(x,y)| c_1 c_2 t^{1/4} \e^{-
{s^2(x,y)}/({2t})},
\end{eqnarray*}
where the last inequality follows from the Cauchy--Schwarz inequality,
Proposition~\ref{prop:bound2} and assumption
(\ref{big assumption}). We thus have
\[
|\rho_2(x,y,t)-\rho_1(x,y,t)|\leq\frac
{c_3}{2}|LC_0(x,y)|t^{1/4}\e^{-{s^2(x,y)}/({2t})}.
\]
Next, assume the induction statement is true for $2,3,\ldots,j-1$.
Then
\begin{eqnarray*}
&&
|\rho_{j+1}(x,y,t)-\rho_{j}(x,y,t)|\\
&&\qquad\leq
\int_0^t\int_{-\infty}^\infty|\kappa^*(x,z,t-\tau)||\rho
_j(z,y,\tau)-\rho_{j-1}(z,y,\tau)|\,dz\,d\tau\\
&&\qquad\leq\int_0^t\int_{-\infty}^\infty\frac{\e^{-
{s^2(x,z)}/({2(t-\tau)})}}{\sqrt{t-\tau}}|LC_0(x,z)|
\frac{c_3}{2^{j-1}}|LC_0(z,y)|\\
&&\qquad\quad\hspace*{34.3pt}{}\times\tau^{1/4}\e^{-
{s^2(z,y)}/({2\tau})} \,dz\,d\tau\\
&&\qquad\leq\frac{c_3}{2^{j-1}}|LC_0(x,y)|\int_0^t\int_{-\infty}^\infty
\frac{\e^{-{s^2(x,z)}/({2(t-\tau)})}}{\sqrt{t-\tau}}\frac{\e
^{-{s^2(z,y)}/({2\tau})}}{\sqrt{\sigma(z)\tau}}\\
&&\qquad\quad\hspace*{112.1pt}{}\times\tau^{3/4}
\frac{|Lq(z)|}{\sqrt{2\pi}q(z)} \,dz\,d\tau\\
&&\qquad\leq\frac{c_3}{2^{j-1}}|LC_0(x,y)| t^{1/4} \e^{-
{s^2(x,y)}/({2t})}t_0^{3/4} \frac{c_1 c_2}{\sqrt{2\pi}}.
\end{eqnarray*}
The last line follows from the
Cauchy--Schwarz inequality and the fact that $\tau^{3/4}\leq
t^{3/4}\leq t_0^{3/4}$. Since
$t_0^{3/4} \frac{c_1 c_2}{\sqrt{2\pi}}\leq\frac{1}{2}$, we obtain
\[
|\rho_{j+1}(x,y,t)-\rho_{j}(x,y,t)|\leq\frac
{c_3}{2^j}|LC_0(x,y)|t^{1/4}\e^{-{s^2(x,y)}/({2t})}.
\]
This establishes (\ref{indo}). Next, we have the bound for all $j\geq1$:
%
%
\begin{eqnarray}\label{thebound}
|\rho_j(x,y,t)|&\leq&|\rho_1(x,y,t)|+\sum_{j=1}^\infty\frac
{c_3}{2^j}|LC_0(x,y)|t^{1/4}\e^{-{s^2(x,y)}/({2t})}\nonumber\\
&\leq&|LC_0(x,y)| \biggl(\frac{1}{\sqrt{t}}+c_3t^{1/4} \biggr)
\e^{-{s^2(x,y)}/({2t})}\\
&\leq&|LC_0(x,y)|\frac{2}{\sqrt{t}}\e^{-
{s^2(x,y)}/({2t})}.\nonumber
\end{eqnarray}
%
In the light of (\ref{thebound}) and (\ref{indo}), the
pointwise limit
\[
\rho(x,y,t)=\lim_{j\rightarrow\infty}\rho_j(x,y,t)
\]
exists on $\mathbb{R}\times\mathbb{R}\times(0,t_0)$.
In addition, $\rho(x,y,t)$ satisfies the limiting equation
\[
0=\kappa^*(x,y,t)+\rho(x,y,t)+\int_0^t\int_{-\infty}^\infty
\kappa^*(x,z,t-\tau)\rho(z,y,\tau)\,dz\,d\tau,
\]
and indeed
%
%
\begin{equation}\label{likeCHK}
\kappa(x,y;t)-\widetilde\kappa(x,y;t)=
\int_0^t \int_{-\infty}^{\infty}\widetilde\kappa(x,z,t-\tau)\rho
(z,y,\tau) \,dz\,d\tau.
\end{equation}
In order to see this, we can apply directly the arguments of Section 5
of \cite{Cohen}
in the case $N=0$; see also Section 1.3 of \cite{Friedman}.
Hence, we can take the limit in (\ref{thebound}) to conclude
%
%
\begin{equation}\label{thelimitbound}
|\rho(x,y,t)|
\leq2 |LC_0(x,y)| t^{-1/2}\e^{-{s^2(x,y)}/({2t})}
\end{equation}
for $t\in(0, t_0]$. The claim of the lemma then follows from
\begin{eqnarray*}
&&
|\kappa(x,y;t)-\widetilde\kappa(x,y;t)|\\
&&\qquad\leq
\int_0^t \int_{-\infty}^{\infty}\widetilde\kappa(x,z,t-\tau
)|\rho(z,y,\tau)| \,dz\,d\tau\\
&&\qquad\leq2\int_0^t\int_{-\infty}^\infty\frac{\e^{-
{s^2(x,z)}/({2(t-\tau)})}}{\sqrt{t-\tau}}C_0(x,z)
\frac{\e^{-{s^2(z,y)}/({2\tau})}}{\sqrt{\tau}} |LC_0(z,y)| \,dz\,d\tau\\
&&\qquad\leq\frac{2}{\sqrt{2\pi}}C_0(x,y)\int_0^t\int_{-\infty}^\infty
\frac{\e^{-{s^2(x,z)}/({2(t-\tau)})}}{\sqrt{t-\tau}}\frac
{\e^{-{s^2(z,y)}/({2\tau})}}{\sqrt{\sigma(z)\tau}}
\frac{|Lq(z)|}{q(z)} \,dz\,d\tau\\
&&\qquad\leq2C_0(x,y)t^{1/4}\e^{-{s^2(x,y)}/({2t})} \frac{c_1c_2}{\sqrt
{2\pi}}=c_3C_0(x,y)t^{1/4}\e^{-{s^2(x,y)}/({2t})}.
\end{eqnarray*}
\upqed\end{pf}
%

%
%

\vspace*{-14pt}

\section{\texorpdfstring{Proof of Theorem \protect\lowercase{\ref{diffusion
theorem}}}{Proof of Theorem 1}}
\label{appendix:C}

Note that (\ref{theorem:bias})
is given by $\int_{-\infty}^\infty\kappa(x,y;t)f(y)\,dy-f(x)$, and
from (\ref{BE}) we have
\begin{eqnarray*}
\frac{\partial}{\partial t}g(x;t)&=&
\int_\mathscr{X} f(y)L^*\kappa(x,y;t) \,dy \\
&=&-\frac{1}{2}\,\frac{d}{dy} \biggl(\frac{f(y)}{p(y)} \biggr)a(y)\kappa(x,y;t)
\Big|_{y\in\partial\mathscr{X}}+ \int_\mathscr{X} \kappa(y,x;t)Lf(x)\,dx.
\end{eqnarray*}
Given that $\mathscr{X}\equiv\mathbb{R}$, Lemma \ref{lem:kernel} gives
$\kappa(x,y;t) |_{y\in\partial\mathscr{X}}\sim\widetilde\kappa
(x,y;t) |_{y=-\infty}^{y=\infty} , t\downarrow0$.
The last term is zero since for fixed $x$,
\[
\lim_{y\rightarrow\pm\infty} \Biggl[\int_y^x\sqrt{\frac{p(s)}{a(s)}}
\,ds \Biggr]^2=\infty,
\]
and hence
$ \lim_{y\rightarrow\pm\infty}\widetilde\kappa(x,y;t)=0 $. We have
\[
g(x;t)=g(x;0)+t\,\frac{\partial}{\partial t}g(x;t) \Big|_{t=0}+O(t^2),
\]
because $g(x;t),t>0$ is smooth (see, e.g., Theorem $\mathrm{IV}\cdot10\cdot1$
in \cite{pde}).
Therefore,
\[
g(x;t)=f(x)+t Lf(x)+O(t^2),
\]
and (\ref{theorem:bias}) and (\ref{diffusion_bias}) follow.
We now proceed to demonstrate (\ref{theorem:var}). First, the second
moment has the behavior
\begin{eqnarray*}
&&
\mathbb E_f[\kappa^2(x,Y;t)]\\
&&\qquad= \int_{\mathscr{X}}f(y)\kappa
^2(x,y;t)\,dy \sim\int_{\mathscr{X}}f(y)\widetilde\kappa^2(x,y;t)\,dy \\
&&\qquad\sim\frac{p^2(x)}{2\pi t\sqrt{p(x)a(x)}}\int_{-\infty}^\infty
\frac{f(y)}{\sqrt{p(y)a(y)}} \e^{-{1/2} [\sqrt{
{2/t}}\int_x^y\sqrt{{p(s)}/{a(s)}} \,ds ]^2}\,dy.
\end{eqnarray*}

We can simplify the last expression by the change of variable
$u=\sqrt{\frac{2}{t}}\times\break\int_x^y\sqrt{\frac{p(s)}{a(s)}} \,ds$. This gives
\[
\frac{p^2(x)}{2\pi\sqrt{2t}\sqrt{p(x)a(x)}}\int_{-\infty}^\infty
\frac{f(y(u,t))}{p(y(u,t))} \e^{-{u^2}/{2}}\,du,
\]
where $y(u,t)=y(u,0)+\sqrt{t}\frac{\partial y}{\partial\sqrt{t}}
|_{t=0}+O(t)=x+u\sqrt{\frac{ta(x)}{2p(x)}}+O(t)$ is a Taylor
expansion of $y(u,t)$ at $\sqrt{t}=0$. Therefore, $\frac
{f(y(u,t))}{p(y(u,t))}\sim\frac{f(x)}{p(x)}$ as $t\downarrow0$, and
\[
\frac{p^2(x)}{2\pi\sqrt{2t}\sqrt{p(x)a(x)}}\int_{-\infty}^\infty
\frac{f(y(u,t))}{p(y(u,t))} \e^{-{u^2}/{2}}\,du\sim\frac
{1}{2\sqrt{\pi t}}f(x)\sqrt{\frac{p(x)}{a(x)}},\qquad t\downarrow0.
\]
Hence, from (\ref{diff est}) we have
\begin{eqnarray*}
\operatorname{Var}_f[g(x;t)]&=&\frac{1}{N} \mathbb E_f[\kappa
^2(x,Y;t)]-\frac
{1}{N}\mathbb E_f[\kappa(x,Y;t)]^2
\\
&\sim&\frac{f(x)}{2N\sqrt{\pi t }
\sigma
(x)},\qquad t\downarrow0,
\end{eqnarray*}
from which (\ref{diffusion_var}) and (\ref{theorem:var}) follow.

\section{Consistency at boundary}
\label{appendix:boundary}

As in \cite{WandJones}, we consider the case where the
support of $f$ is $[0,\infty]$. The
consistency of the estimator near $x=0$ is analyzed by considering the
pointwise bias of estimator (\ref{diff est}) at a
point $x_N$ such that $x_N$ is
$O(\sqrt{t_N})$ away from the boundary, that is, $x_N$ is
approaching the boundary at the same rate at which the bandwidth is
approaching $0$.
We then have the following result, which shows that the diffusion
estimator (\ref{diff est}), and hence its special case (\ref{theta
estimator}), is consistent at the boundaries.
\begin{proposition}
Let $\mathscr{X}\equiv[0,\infty]$, and assume that $x=x_N=\alpha
\sqrt{t_N}$ for some constant $\alpha\in[0, 1]$, where $\lim
_{N\rightarrow\infty}t_N=0$ and $\lim_{N\rightarrow\infty} N\sqrt
{t_N}=\infty$.
Then for the diffusion estimator (\ref{diff est}) we have
\[
\mathbb E_f g(x_N;t)=f(x_N)+O\bigl(\sqrt{t_N}\bigr),\qquad N\rightarrow\infty.
\]
Hence, the diffusion estimator (\ref{diff est}) is consistent at the
boundaries.
\end{proposition}
\begin{pf}
First, we differentiate both sides of $ \mathbb E_f g(x;t)=\int_0^1
f(y)\kappa(x;y;t) \,dy$
with respect to $t$ and use (\ref{BE}) to obtain
\begin{eqnarray*}
\frac{\partial}{\partial t} \mathbb E_f g(x;t)&=&\int_0^\infty
f(y)\,\frac{\partial}{\partial t}\kappa(x;y;t) \,dy\\
&=& \int_0^\infty
f(y)L^*\kappa(x;y;t)\,dy\\
&=& - \frac{1}{2} \biggl(\frac{f(y)}{p(y)} \biggr)'
a(y) \kappa(x;y;t) \Big|_{y=0}^{y=\infty} +\int_0^\infty\kappa
(x;y;t)Lf(y) \,d y.
\end{eqnarray*}
Second, we show that
$\kappa(\alpha\sqrt{t_N};0;t_N)=O(t^{-1/2})$ and
$\lim_{y\rightarrow\infty}\kappa(\alpha\sqrt{t_N};y;t_N)=o(1) $,
and $\int_0^1 \kappa(x;y;t_N)Lf(y) \,d y=O(1)$ as $N\rightarrow\infty$. To this end, we
consider the small bandwidth behavior of $\kappa$. It is
easy to verify using Lemma \ref{lem:kernel} that the \textit{boundary
kernel}
\[
\kappa_B(x,y;t)=\widetilde\kappa(x,y;t)+\widetilde\kappa(x,-y;t)
\]
satisfies
\[
\frac{\partial}{\partial t}\kappa_B(x,y;t)=L^*\kappa_B(x,y;t)+O \bigl(\e
^{-{s^2(x,y)}/({2t})} t^{-1/2} \bigr),\qquad t\downarrow0,
\]
on
$x,y\in\mathbb{R}$ with initial condition $\kappa_B(x,y;0)=\delta(x-y)$.
In addition, the boundary kernel satisfies the condition
$\frac{\partial}{\partial y}\kappa_B(x,y;t) |_{y=0}=0$, and
therefore $\kappa_B$ describes the small bandwidth
asymptotics of the solution of the PDE
(\ref{BE}) on the domain $x,y \in[0,\infty)$ with boundary
condition $\frac{\partial}{\partial
y}\kappa(x,y;t) |_{y=0}=0$.
Hence, we have
\[
\kappa\bigl(\alpha\sqrt{t};0;t\bigr)\sim\kappa_B\bigl(\alpha\sqrt
{t};0;t\bigr)=\operatorname{const.} t^{-1/2}\e^{O(\sqrt{t})},\qquad t\downarrow0,
\]
and
\[
\lim_{y\rightarrow\infty}\kappa_B\bigl(\alpha\sqrt{t};y;t\bigr)=0,\qquad t>0.
\]
Therefore,
\[
\frac{\partial}{\partial t} \mathbb E_f g(x_N;t_N) =
o(1)-O(t_N^{-1/2}) ,\qquad
N\rightarrow\infty,
\]
or
\[
\frac{\mathbb E_fg(x_N;t_N)-\mathbb E_fg(x_N;0)}{t_N}+O(t_N)=O(t_N^{-1/2})+O(1)
,\qquad N \rightarrow\infty,
\]
which, after rearranging, gives
\[
\mathbb E_fg(x_N;t_N)=f(x_N)+O\bigl(\sqrt{t_N}\bigr),\qquad N\rightarrow\infty.
\]
\upqed\end{pf}

\vspace*{-14pt}

\section{Bandwidth selection in higher dimensions}
\label{appendix:E}

Algorithm \ref{my plug-in} can be extended to two dimensions for the
estimation of a p.d.f. $f(\mathbf{x})$ on $\mathbb{R}^2$. Assuming
a Gaussian kernel
\[
\phi(\mathbf{x},\mathbf{y};t)=\frac{1}{2\pi t} \e^{-
{(\mathbf
{x}-\mathbf{y})^T(\mathbf{x}-\mathbf{y})}/({2t})},
\]
where $\mathbf{x}=[x_1,x_2]^T$ and $\mathbf{y}=[y_1,y_2]^T$,
the asymptotically optimal squared bandwidth is given by
(\cite{WandJones}, page 99)
%
\[
t^*= \bigl(2\pi N (\psi_{0,2}+ \psi_{2,0}+2\psi_{1,1} ) \bigr)^{-1/3},
\]
where
%
%
\begin{eqnarray}
\label{2dim functional}
\psi_{i,j}&=& (-1)^{i+j}\int_{\mathbb{R}^2} f(\mathbf{x})\, \frac
{\partial^{2(i+j)}}{\partial x_1^{2i}\,\partial x_2^{2j}}f(\mathbf{x})\,
d\mathbf{x},\qquad i,j\in\mathbb{N}^+,\nonumber\\[-8pt]\\[-8pt]
&=&\int\biggl(\frac{\partial^{(i+j)} }{\partial x_1^i\,\partial x_2^j}
f(\mathbf{x}) \biggr)^2\,d\mathbf{x}.\nonumber
\end{eqnarray}
Note that our definition of $\psi$ differs slightly from the
definition of $\psi$ in \cite{WandJones}. Here the partial
derivatives under the integral sign are applied $2(i+j)$ times, while
in \cite{WandJones} they are applied $(i+j)$ times.
Similar to the one-dimensional case,
there are two viable plug-in estimators for $\psi_{i,j}$. The first
one is derived from the first line of (\ref{2dim functional}):
%
%
\begin{equation}
\widetilde\psi_{i,j}=\frac{(-1)^{i+j}}{N^2}\sum_{k=1}^N\sum
_{m=1}^N\frac{\partial^{2(i+j)}}{\partial x_1^{2i}\,\partial
x_2^{2j}}\phi(\mathbf{X}_m,\mathbf{X}_k;t_{i,j}),
\end{equation}
and the second one is derived from the second line of
(\ref{2dim functional}):
%
%
\begin{eqnarray}
\label{2dim est}
\widehat\psi_{i,j}&=&\frac{1}{N^2}\sum_{k=1}^N\sum_{m=1}^N\int
\frac{\partial^{(i+j)}}{\partial x_1^{i}\,\partial x_2^{j}}\phi
(\mathbf{x},\mathbf{X}_m;t_{i,j})\,\frac{\partial^{(i+j)}}{\partial
x_1^{i}\,\partial x_2^{j}}\phi(\mathbf{x},\mathbf{X}_k;t_{i,j})
\,d\mathbf{x}\nonumber\\[-8pt]\\[-8pt]
&=&\frac{(-1)^{i+j}}{N^2}\sum_{k=1}^N\sum_{m=1}^N\frac{\partial
^{2(i+j)}}{\partial x_1^{2i}\,\partial
x_2^{2j}}\phi(\mathbf{X}_m,\mathbf{X}_k;2 t_{i,j}).\nonumber
\end{eqnarray}
The asymptotic expansion of the squared bias of estimator
$\widetilde\psi_{i,j}$ is given by (\cite{WandJones}, page 113)
%
%
\begin{eqnarray}
\label{squared bias}
&&
(\mathbb E_f[\widetilde\psi_{i,j}]-\psi_{i,j} )^2\nonumber\\[-8pt]\\[-8pt]
&&\qquad\sim\biggl(\frac{q(i)
q(j)}{N t^{i+j+1}_{i,j}}+\frac{t_{i,j}}{2} (\psi_{i+1,j}+\psi
_{i,j+1} ) \biggr)^2,\qquad N\rightarrow\infty,\nonumber
\end{eqnarray}
where
\[
q(j)=\cases{(-1)^j \dfrac{1\times3\times5\times\cdots\times
(2j-1)}{\sqrt{2\pi}}, &\quad $j\geq
1$,\cr
\dfrac{1}{\sqrt{2\pi}}, &\quad $j=0$.}
\]
Thus, we have
%
%
\begin{eqnarray}
\label{squared bias 2}
&&(\mathbb E_f[\widehat\psi_{i,j}]-\psi_{i,j} )^2\nonumber\\[-8pt]\\[-8pt]
&&\qquad\sim\biggl(\frac{q(i) q(j)}{N
(2 t_{i,j})^{i+j+1}}+t_{i,j} (\psi_{i+1,j}+\psi_{i,j+1} )
\biggr)^2,\qquad
N\rightarrow\infty.\nonumber
\end{eqnarray}
For both estimators the squared bias is the dominant term
in the asymptotic mean squared error, because the variance is of the
order $O(N^{-2}t^{-2i-2j-1})$.
It follows that both estimators will have the same leading asymptotic
mean square error term provided that
%
%
\begin{equation}
\label{2 dim t*} t_{i,j}= \biggl(\frac{1+2^{-i-j-1}}{3} \frac{-2 q(i)
q(j)}{N(\psi_{i+1,j}+\psi_{i,j+1})} \biggr)^{1/(2+i+j)}.
\end{equation}
We estimate $t_{i,j}$ via
%
%
\begin{equation}
\label{2 dim t* est} \hat t_{i,j}= \biggl(\frac{1+2^{-i-j-1}}{3} \frac{-2
q(i) q(j)}{N(\widehat{\psi}_{i+1,j}+\widehat{\psi}_{i,j+1})} \biggr)^{1/(2+i+j)}.
\end{equation}
%
Thus, estimation of $\psi_{i,j}$ requires
estimation of $\psi_{i,j+1}$ and $\psi_{i+1,j}$, which in turn
requires estimation of $\psi_{i+2,j}, \psi_{i+1,j+1}, \psi_{i,j+2}$
and so on applying formula (\ref{2 dim t* est}), recursively. Observe
that to estimate all
$\psi_{i,j}$ for which $i+j=k$, that is,
$\{\psi_{i,j}\dvtx i+j=k\}$, we need estimates of
all $\{\psi_{i,j}\dvtx i+j=k+1\}$. For example, from
formula (\ref{2 dim t* est}) we can see that estimation of $t_{2,0},
t_{1,1}, t_{0,2}$ requires estimation of $t_{3,0}, t_{2,1}, t_{1,2}, t_{0,3}$.

For a given integer $k\geq3$, we define the function $\gamma(t)$ as
follows. Given an input $t>0$:
\begin{enumerate}
\item Set $\hat t_{i,j}=t$ for all $i+j=k$.
\item Use the set $\{\hat t_{i,j}\dvtx i+j=k\}$ to compute all functionals
$\{\widehat\psi_{i,j}\dvtx i+j=k\}$ via
(\ref{2dim est}).
\item Use $\{\widehat\psi_{i,j}\dvtx i+j=k\}$ to compute
$\{\hat t_{i,j}\dvtx i+j=k-1\}$ via (\ref{2 dim t* est}).
%
\item If $k=2$ go to step 5; otherwise set $k:=k-1$ and repeat from
step 2.
\item Use $\{\widehat\psi_{i,j}\dvtx i+j=2\}$
to output
\[
\gamma(t)= \bigl(2\pi N (\widehat\psi_{0,2}+ \widehat\psi
_{2,0}+2\widehat\psi_{1,1} ) \bigr)^{-1/3}.
\]
\end{enumerate}
The bandwidth selection rule simply consists of solving
the equation $\gamma(t)=t$ for a given $k\geq3$ via either
the fixed point iteration in Algorithm \ref{my plug-in}
(ignoring step 4) or by using Newton's method. We obtain excellent numerical
results for $k=4$ or $k=5$. Higher values of $k$ did
not change the value of $t$ in any significant way, but only increased
the computational
cost of evaluating the function $\gamma(t)$. Again
note that this appears to be the first successful
plug-in bandwidth selection rule that does not involve
any arbitrary reference rules, but it is purely data-driven.
An efficient Matlab implementation of
the bandwidth selection rule described
here, and using the two-dimensional discrete cosine
transform, can be downloaded freely from~\cite{Matlab}. The
Matlab implementation takes an additional step in
which, once a fixed point of $\gamma(t)$ has been found, the final set
of estimates
$\{\widehat\psi_{i,j}\dvtx i+j=2\}$ is used to compute the entries
$\sqrt
{t_{X_1}}$ and $\sqrt{t_{X_2}}$ of the optimal diagonal bandwidth
matrix (\cite{WandJones}, page~111) for a Gaussian kernel of the form
\[
\frac{1}{2\pi\sqrt{t_{X_1} t_{X_2}}} \e^{-
{(x_1-y_1)^2}/({2t_{X_1}})-{(x_2-y_2)^2}/({2t_{X_2}})}.
\]
These entries are estimated via the formulas
\[
t_{X_1}= \biggl(\frac{\widehat\psi_{0,2}^{3/4}}{4\pi N \widehat\psi
_{2,0}^{3/4} (\widehat\psi_{1,1}+\sqrt{\widehat\psi_{2,0} \widehat
\psi_{0,2}} )} \biggr)^{1/3}
\]
and
\[
t_{X_2}= \biggl(\frac{\widehat\psi_{2,0}^{3/4}}{4\pi N \widehat\psi
_{0,2}^{3/4} (\widehat\psi_{1,1}+\sqrt{\widehat\psi_{2,0} \widehat
\psi_{0,2}} )} \biggr)^{1/3}.
\]
\end{appendix}

\printaddresses

\end{document}